\documentclass[reqno,12pt,a4paper]{amsart}

\usepackage{mathrsfs}
\usepackage{amssymb}
\usepackage{amsfonts}
\usepackage{latexsym}
\usepackage{amsthm}
\usepackage{graphicx}

\def\lb{\label}

\newcommand{\er}[1]{\textrm{(\ref{#1})}}

\begin{document}

%%%%%%%%%% Some definitions %%%%%%%%%%

%%%%%%%% Equations, theorems %%%%%%%%%
\renewcommand{\theequation}{\arabic{section}.\arabic{equation}}
\theoremstyle{plain}
\newtheorem{theorem}{\bf Theorem}[section]
\newtheorem{lemma}[theorem]{\bf Lemma}
\newtheorem{corollary}[theorem]{\bf Corollary}
\newtheorem{proposition}[theorem]{\bf Proposition}
\newtheorem{definition}[theorem]{\bf Definition}

\newtheorem{remark}[theorem]{\bf Remark}

%%%%% Alphabet %%%%%
\def\a{\alpha}  \def\cA{{\mathcal A}}     \def\bA{{\bf A}}  \def\mA{{\mathscr A}}
\def\b{\beta}   \def\cB{{\mathcal B}}     \def\bB{{\bf B}}  \def\mB{{\mathscr B}}
\def\g{\gamma}  \def\cC{{\mathcal C}}     \def\bC{{\bf C}}  \def\mC{{\mathscr C}}
\def\G{\Gamma}  \def\cD{{\mathcal D}}     \def\bD{{\bf D}}  \def\mD{{\mathscr D}}
\def\d{\delta}  \def\cE{{\mathcal E}}     \def\bE{{\bf E}}  \def\mE{{\mathscr E}}
\def\D{\Delta}  \def\cF{{\mathcal F}}     \def\bF{{\bf F}}  \def\mF{{\mathscr F}}
\def\c{\chi}    \def\cG{{\mathcal G}}     \def\bG{{\bf G}}  \def\mG{{\mathscr G}}
\def\z{\zeta}   \def\cH{{\mathcal H}}     \def\bH{{\bf H}}  \def\mH{{\mathscr H}}
\def\e{\eta}    \def\cI{{\mathcal I}}     \def\bI{{\bf I}}  \def\mI{{\mathscr I}}
\def\p{\psi}    \def\cJ{{\mathcal J}}     \def\bJ{{\bf J}}  \def\mJ{{\mathscr J}}
\def\vT{\Theta} \def\cK{{\mathcal K}}     \def\bK{{\bf K}}  \def\mK{{\mathscr K}}
\def\k{\kappa}  \def\cL{{\mathcal L}}     \def\bL{{\bf L}}  \def\mL{{\mathscr L}}
\def\l{\lambda} \def\cM{{\mathcal M}}     \def\bM{{\bf M}}  \def\mM{{\mathscr M}}
\def\L{\Lambda} \def\cN{{\mathcal N}}     \def\bN{{\bf N}}  \def\mN{{\mathscr N}}
\def\m{\mu}     \def\cO{{\mathcal O}}     \def\bO{{\bf O}}  \def\mO{{\mathscr O}}
\def\n{\nu}     \def\cP{{\mathcal P}}     \def\bP{{\bf P}}  \def\mP{{\mathscr P}}
\def\r{\rho}    \def\cQ{{\mathcal Q}}     \def\bQ{{\bf Q}}  \def\mQ{{\mathscr Q}}
\def\s{\sigma}  \def\cR{{\mathcal R}}     \def\bR{{\bf R}}  \def\mR{{\mathscr R}}
\def\S{\Sigma}  \def\cS{{\mathcal S}}     \def\bS{{\bf S}}  \def\mS{{\mathscr S}}
\def\t{\tau}    \def\cT{{\mathcal T}}     \def\bT{{\bf T}}  \def\mT{{\mathscr T}}
\def\f{\phi}    \def\cU{{\mathcal U}}     \def\bU{{\bf U}}  \def\mU{{\mathscr U}}
\def\F{\Phi}    \def\cV{{\mathcal V}}     \def\bV{{\bf V}}  \def\mV{{\mathscr V}}
\def\P{\Psi}    \def\cW{{\mathcal W}}     \def\bW{{\bf W}}  \def\mW{{\mathscr W}}
\def\o{\omega}  \def\cX{{\mathcal X}}     \def\bX{{\bf X}}  \def\mX{{\mathscr X}}
\def\x{\xi}     \def\cY{{\mathcal Y}}     \def\bY{{\bf Y}}  \def\mY{{\mathscr Y}}
\def\X{\Xi}     \def\cZ{{\mathcal Z}}     \def\bZ{{\bf Z}}  \def\mZ{{\mathscr Z}}
%*********************
\def\be{{\bf e}}
\def\bv{{\bf v}} \def\bu{{\bf u}}
\def\Om{\Omega}
%************************
\def\bbD{\pmb \Delta}
\def\mm{\mathrm m}
\def\mn{\mathrm n}
\def\vr{\varrho}
%*************************

\newcommand{\mc}{\mathscr {c}}

\newcommand{\gA}{\mathfrak{A}}          \newcommand{\ga}{\mathfrak{a}}
\newcommand{\gB}{\mathfrak{B}}          \newcommand{\gb}{\mathfrak{b}}
\newcommand{\gC}{\mathfrak{C}}          \newcommand{\gc}{\mathfrak{c}}
\newcommand{\gD}{\mathfrak{D}}          \newcommand{\gd}{\mathfrak{d}}
\newcommand{\gE}{\mathfrak{E}}
\newcommand{\gF}{\mathfrak{F}}           \newcommand{\gf}{\mathfrak{f}}
\newcommand{\gG}{\mathfrak{G}}           %\newcommand{\gg}{\mathfrak{g}}
\newcommand{\gH}{\mathfrak{H}}           \newcommand{\gh}{\mathfrak{h}}
\newcommand{\gI}{\mathfrak{I}}           \newcommand{\gi}{\mathfrak{i}}
\newcommand{\gJ}{\mathfrak{J}}           \newcommand{\gj}{\mathfrak{j}}
\newcommand{\gK}{\mathfrak{K}}            \newcommand{\gk}{\mathfrak{k}}
\newcommand{\gL}{\mathfrak{L}}            \newcommand{\gl}{\mathfrak{l}}
\newcommand{\gM}{\mathfrak{M}}            \newcommand{\gm}{\mathfrak{m}}
\newcommand{\gN}{\mathfrak{N}}            \newcommand{\gn}{\mathfrak{n}}
\newcommand{\gO}{\mathfrak{O}}
\newcommand{\gP}{\mathfrak{P}}             \newcommand{\gp}{\mathfrak{p}}
\newcommand{\gQ}{\mathfrak{Q}}             \newcommand{\gq}{\mathfrak{q}}
\newcommand{\gR}{\mathfrak{R}}             \newcommand{\gr}{\mathfrak{r}}
\newcommand{\gS}{\mathfrak{S}}              \newcommand{\gs}{\mathfrak{s}}
\newcommand{\gT}{\mathfrak{T}}             \newcommand{\gt}{\mathfrak{t}}
\newcommand{\gU}{\mathfrak{U}}             \newcommand{\gu}{\mathfrak{u}}
\newcommand{\gV}{\mathfrak{V}}             \newcommand{\gv}{\mathfrak{v}}
\newcommand{\gW}{\mathfrak{W}}             \newcommand{\gw}{\mathfrak{w}}
\newcommand{\gX}{\mathfrak{X}}               \newcommand{\gx}{\mathfrak{x}}
\newcommand{\gY}{\mathfrak{Y}}              \newcommand{\gy}{\mathfrak{y}}
\newcommand{\gZ}{\mathfrak{Z}}             \newcommand{\gz}{\mathfrak{z}}

\def\ve{\varepsilon}   \def\vt{\vartheta}    \def\vp{\varphi}    \def\vk{\varkappa}

\def\A{{\mathbb A}} \def\B{{\mathbb B}} \def\C{{\mathbb C}}
\def\dD{{\mathbb D}} \def\E{{\mathbb E}} \def\dF{{\mathbb F}} \def\dG{{\mathbb G}} \def\H{{\mathbb H}}\def\I{{\mathbb I}} \def\J{{\mathbb J}} \def\K{{\mathbb K}} \def\dL{{\mathbb L}}\def\M{{\mathbb M}} \def\N{{\mathbb N}} \def\O{{\mathbb O}} \def\dP{{\mathbb P}} \def\R{{\mathbb R}}\def\S{{\mathbb S}} \def\T{{\mathbb T}} \def\U{{\mathbb U}} \def\V{{\mathbb V}}\def\W{{\mathbb W}} \def\X{{\mathbb X}} \def\Y{{\mathbb Y}} \def\Z{{\mathbb Z}}

%%%%% Arrows %%%%%

\def\la{\leftarrow}              \def\ra{\rightarrow}            \def\Ra{\Rightarrow}
\def\ua{\uparrow}                \def\da{\downarrow}
\def\lra{\leftrightarrow}        \def\Lra{\Leftrightarrow}

%%%%% Typography %%%%%

\def\lt{\biggl}                  \def\rt{\biggr}
\def\ol{\overline}               \def\wt{\widetilde}
\def\ul{\underline}
\def\no{\noindent}

%%%%% Math signs %%%%%

\let\ge\geqslant                 \let\le\leqslant
\def\lan{\langle}                \def\ran{\rangle}
\def\/{\over}                    \def\iy{\infty}
\def\sm{\setminus}               \def\es{\emptyset}
\def\ss{\subset}                 \def\ts{\times}
\def\pa{\partial}                \def\os{\oplus}
\def\om{\ominus}                 \def\ev{\equiv}
\def\iint{\int\!\!\!\int}        \def\iintt{\mathop{\int\!\!\int\!\!\dots\!\!\int}\limits}
\def\el2{\ell^{\,2}}             \def\1{1\!\!1}
\def\sh{\sharp}
\def\wh{\widehat}
\def\bs{\backslash}
\def\intl{\int\limits}
%%%%% Math operations %%%%%

\def\na{\mathop{\mathrm{\nabla}}\nolimits}
\def\sh{\mathop{\mathrm{sh}}\nolimits}
\def\ch{\mathop{\mathrm{ch}}\nolimits}
\def\where{\mathop{\mathrm{where}}\nolimits}
\def\all{\mathop{\mathrm{all}}\nolimits}
\def\as{\mathop{\mathrm{as}}\nolimits}
\def\Area{\mathop{\mathrm{Area}}\nolimits}
\def\arg{\mathop{\mathrm{arg}}\nolimits}
\def\const{\mathop{\mathrm{const}}\nolimits}
\def\det{\mathop{\mathrm{det}}\nolimits}
\def\diag{\mathop{\mathrm{diag}}\nolimits}
\def\diam{\mathop{\mathrm{diam}}\nolimits}
\def\dim{\mathop{\mathrm{dim}}\nolimits}
\def\dist{\mathop{\mathrm{dist}}\nolimits}
\def\Im{\mathop{\mathrm{Im}}\nolimits}
\def\Iso{\mathop{\mathrm{Iso}}\nolimits}
\def\Ker{\mathop{\mathrm{Ker}}\nolimits}
\def\Lip{\mathop{\mathrm{Lip}}\nolimits}
\def\rank{\mathop{\mathrm{rank}}\limits}
\def\Ran{\mathop{\mathrm{Ran}}\nolimits}
\def\Re{\mathop{\mathrm{Re}}\nolimits}
\def\Res{\mathop{\mathrm{Res}}\nolimits}
\def\res{\mathop{\mathrm{res}}\limits}
\def\sign{\mathop{\mathrm{sign}}\nolimits}
\def\span{\mathop{\mathrm{span}}\nolimits}
\def\supp{\mathop{\mathrm{supp}}\nolimits}
\def\Tr{\mathop{\mathrm{Tr}}\nolimits}
\def\BBox{\hspace{1mm}\vrule height6pt width5.5pt depth0pt \hspace{6pt}}

%%%%%%%%%%%%% specialities %%%%%%%%%%%%%%

\newcommand\nh[2]{\widehat{#1}\vphantom{#1}^{(#2)}}
%{{\mathop{#1}\limits^\wedge}\vphantom{#1}^{(#2)}}
\def\dia{\diamond}

\def\Oplus{\bigoplus\nolimits}

%%%%%%%%%%% End of definitions %%%%%%%%%%

%%%%% OLD OLD OLD

\def\qqq{\qquad}
\def\qq{\quad}
\let\ge\geqslant
\let\le\leqslant
\let\geq\geqslant
\let\leq\leqslant
\newcommand{\ca}{\begin{cases}}
\newcommand{\ac}{\end{cases}}
\newcommand{\ma}{\begin{pmatrix}}
\newcommand{\am}{\end{pmatrix}}
\renewcommand{\[}{\begin{equation}}
\renewcommand{\]}{\end{equation}}
\def\eq{\begin{equation}}
\def\qe{\end{equation}}
\def\[{\begin{equation}}
\def\bu{\bullet}

\title[{Laplacians on periodic graphs with guides}]
{Laplacians on periodic graphs with guides}

\date{\today}
\author[Evgeny Korotyaev]{Evgeny Korotyaev}
\address{Saint-Petersburg State University, Universitetskaya nab. 7/9,
 St. Petersburg, 199034, Russia,
\ korotyaev@gmail.com, \
e.korotyaev@spbu.ru,}
\author[Natalia Saburova]{Natalia Saburova}
\address{Northern (Arctic) Federal University, Severnaya Dvina Emb. 17,
Arkhangelsk, 163002, Russia,
 \ n.saburova@gmail.com, \ n.saburova@narfu.ru}


\begin{abstract}
We consider Laplace operators on periodic discrete graphs perturbed
by guides, i.e., graphs which are periodic in some directions and
finite in other ones.  The spectrum of the Laplacian on the
unperturbed graph is a union of a finite number of non-degenerate
bands and eigenvalues of infinite multiplicity. We show that the
spectrum of the perturbed  Laplacian consists of the unperturbed one
plus the additional so-called guided spectrum which is a union of a
finite number of bands.  We estimate the position of the guided
bands and their length in terms of geometric parameters of the
graph. We also determine the asymptotics of the guided bands for
guides with large multiplicity of edges. Moreover, we show that the
possible number of guided bands, their length and position can be
rather arbitrary for some specific periodic graphs with guides.
\end{abstract}

\maketitle

\vskip 0.25cm

\section {\lb{Sec1}Introduction}
\setcounter{equation}{0} Laplacians on periodic discrete graphs have
attracted a lot of attention due to their applications to the study
of electronic properties of real crystalline structures, see, e.g.,
 \cite{H02},  \cite{NG04} and the survey
\cite{CGPNG09}. However, the arrangement of atoms or molecules in
most crystalline materials is not perfect. The regular patterns are
interrupted by crystalline defects. These defects are the most
important features of the engineering material and are manipulated
to control its behavior.

We consider  Laplace operators on periodic discrete graphs perturbed
by guides (i.e., graphs which are periodic in some directions and finite
in others). For example,  a guide is a periodic graph  embedded
into a  strip in the case of planar graphs.  It is well known that
the spectrum of discrete Laplacians on periodic graphs has a band
structure with a finite number of flat bands (eigenvalues of
infinite multiplicity) \cite{HN09}, \cite{HS04}, \cite{KS14},
\cite{RR07}. The spectrum of the Laplacian on the perturbed graph
consists of the spectrum of the Laplacian on the unperturbed
periodic graph plus the so-called \emph{guided} spectrum. The
additional guided spectrum is a union of a finite number of bands
and the corresponding wave-functions are mainly located along the
guides. In our paper we study guided spectra of Laplacians. We
describe our main goals:

%\medskip

\no $\bu $ to estimate the position of guided bands and their
lengths in terms of  geometric parameters of graphs;

\no $\bu $ to determine asymptotics of the guided spectrum for
guides with large multiplicity of edges;

\no $\bu $ to show that a possible number of guided bands
(including flat bands), their length and positions can be rather
arbitrary for some specific periodic graphs with guides.

\subsection{Discrete Laplacians on periodic graphs}
Let $\G=(V,\cE)$ be a connected infinite graph, possibly  having
loops and multiple edges, where $V$ is the set of its vertices and
$\cE$ is the set of its unoriented edges. From the set $\cE$ we
construct  the set $\cA$ of oriented edges by considering each edge
in $\cE$ to have two orientations. An edge starting at a vertex $u$
and ending at a vertex $v$ from $V$ will be denoted as the ordered
pair $(u,v)\in\cA$. Vertices $u,v\in V$ will be called
\emph{adjacent} and denoted by $u\sim v$, if $(u,v)\in \cA$. We
define the degree $\vk_v$ of the vertex $v\in V$ as the number of
all edges from $\cA$ starting at $v$. We consider graphs with uniformly bounded degrees.

Let $\ell^2(V)$ be the
Hilbert space of all functions $f:V\to \C$ equipped
with the norm
$$
\|f\|^2_{\ell^2(V)}=\sum_{v\in V}|f(v)|^2<\infty.
$$
We define the discrete Laplacian (i.e., the combinatorial  Laplace
operator) $\D$ on $\ell^2(V)$ by
\[
\lb{DOLN} \big(\D
f\big)(v)=\sum_{\be=(v,u)\in\cA}\big(f(v)-f(u)\big), \qqq
f=(f(v))_{v\in V}\in\ell^2(V),
\]
where the sum is taken over all oriented edges starting at the
vertex $v\in V$. It is well known, see, e.g., \cite{M91}, that $\D$
is self-adjoint and its spectrum satisfies: \emph{the point 0
belongs to the spectrum $\s(\D)$ containing in $[0,2\vk_+]$, i.e.,}
\[
\lb{bf}
\begin{aligned}
0\in\s(\D)\subset[0,2\vk_+],\qqq
\textrm{where}\qqq
\vk_+=\sup_{v\in V}\vk_v<\infty.
\end{aligned}
\]
We consider a $\Z^{\wt d}$-periodic graph $\G_0=(V_0,\cE_0)$, i.e.,
a graph satisfying the following conditions:
\begin{itemize}
  \item[1)] \emph{$\G_0$ is equipped with an action of the free abelian group $\Z^{\wt d}$;}
  \item[2)] \emph{the quotient graph  $\G_*=(V_*,\cE_*)=\G_0/\Z^{\wt d}$ is finite.}
\end{itemize}

We assume that the graphs are embedded into Euclidean space, since
in many applications such a natural embedding exists. For example,
in the tight-binding approximation real crystalline structures are
modeled as discrete graphs embedded into $\R^d$ ($d=2,3$) and
consisting of vertices (points representing positions of atoms) and
edges (representing chemical bonding of atoms), by ignoring the
physical characters of atoms and bonds that may be different from
one another. But all results of the paper stay valid in the case
of abstract periodic graphs (without the assumption of graph embedding
into Euclidean space).

For a periodic graph $\G_0$ embedded into the space $\R^{\wt d}$,
the quotient graph  $\G_0/\Z^{\wt d}$ is a graph on the $\wt
d$-dimensional  torus $\R^{\wt d}/\Z^{\wt d}$. Due to the
definition, the graph $\G_0$ is invariant under translations through
vectors $a_1,\ldots,a_{\wt d}$\, which generate the group $\Z^{\wt
d}$:
$$
\G_0+a_s=\G_0,\qqq \forall\, s\in\N_{\wt d}\,.
$$
Here and below for each integer
$m$ the set $\N_m$ is given by
\[
\N_{m}=\{1,\ldots,m\,\}.
\]
We will call the vectors $a_1,\ldots,a_{\wt d}$ \emph{the periods of
the graph}  $\G_0$. In the space $\R^{\wt d}$ we consider a
coordinate system with the origin at some point $O$ and with the
basis $a_1,\ldots,a_{\wt d}$. Below the coordinates of all graph
vertices will be expressed  in this coordinate system.

Let $A$ be a self-adjoint operator, we denote by $\sigma(A)$, $\sigma_{ac}(A)$, $\sigma_{p}(A)$, and $\sigma_{fb}(A)$ its spectrum, absolutely continuous spectrum, point spectrum (eigenvalues of finite multiplicity), and the flat band spectrum (eigenvalues of infinite multiplicity), respectively.

We consider the Laplacian defined by \er{DOLN} on the periodic graph
$\G_0$  as an \emph{unperturbed operator} and denote it by $\D_0$.
It is well known that the spectrum $\s(\D_0)$ of the Laplacian on
periodic graphs is a union of $\n$ spectral bands $\s_n(\D_0)$:
\[
\lb{sH0}
\s(\D_0)=\bigcup_{n=1}^\n\s_n(\D_0)=\s_{ac}(\D_0)\cup\s_{fb}(\D_0),
\]
where $\n=\# V_*$ is the number of vertices of the quotient graph
$\G_*$, the absolutely continuous spectrum $\s_{ac}(\D_0)$ consists
of non-degenerate bands $\s_n(\D_0)$. Note that each flat band is a degenerate band. The spectrum $\s(\D_0)$ is a subset of the interval $[0,\vr]$:
\[
\lb{mm} \s(\D_0)\ss [0,\vr],\qqq \inf \s(\D_0)=0,\qqq \vr=\sup
\s(\D_0).
\]
\subsection{Results overview} There are results about spectral  properties
of the Schr\"odinger operator $H_0=\D_0+W$ with a periodic potential  $W$.
The decomposition of the operator $H_0$ into a constant fiber direct
integral was obtained in \cite{HN09}, \cite{HS04}, \cite{RR07}
without an exact form of fiber operators and in \cite{KS14},
\cite{KS17} with an exact form of fiber operators. In particular,
this yields the band-gap structure of the spectrum of the operator
$H_0$. In \cite{GKT93} the authors described different properties of Schr\"odinger
operators with periodic potentials on the lattice $\Z^2$, the
simplest $\Z^2$-periodic graph. In \cite{LP08},
%simplest example of $\Z^2$-periodic graphs. In \cite{LP08},
\cite{KS15} the positions of the spectral bands of the Laplacians
were estimated in terms of eigenvalues of the operator on finite
graphs (the so-called eigenvalue bracketing). The estimate of the
total length of all bands $\s_n(H_0)$ given by
\[
\lb{eq.7}
\sum_{n=1}^{\n}|\s_n(H_0)|\le 2\b,
\]
was obtained in \cite{KS14}; where $\b=\#\cE_*+1-\n$ is the
so-called Betti number, $\#\cE_*$ is the number of edges of the
quotient graph $\G_*$. Moreover, a global variation of the Lebesgue
measure of the spectrum and a global variation of the gap-length in
terms of potentials and geometric parameters of the graph were
determined. Note that the estimate \er{eq.7} also holds true for
magnetic Schr\"odinger operators with periodic magnetic and electric
potentials (see \cite{KS17}). Estimates of the Lebesgue measure of
the spectrum of $H_0$ in terms of eigenvalues of Dirichlet and
Neumann operators on a fundamental domain of the periodic graph were
described in \cite{KS15}. Estimates of effective masses, associated
with the ends of each spectral band, in terms of geometric
parameters of the graphs were obtained in \cite{KS16}. Moreover, in
the case of the bottom of the spectrum two-sided estimates on the
effective mass in terms of geometric parameters of the graphs were
determined. The proof of all these results in
\cite{KS14}-\cite{KS17} is based on Floquet theory and the exact
form of fiber Schr\"odinger operators from \cite{KS14}, \cite{KS17}.  The spectra
of the discrete Schr\"odinger operators on graphene nano-tubes and
nano-ribbons in external fields were discussed in \cite{KK10},
\cite{KK10a}. The spectrum of discrete magnetic Laplacians on some planar graphs (the hexagonal lattice, the kagome lattice and so on) was described in \cite{HKR16} and see  the references therein.

Discrete Laplacians  for some class of periodic graphs with
compact perturbations including the square, triangular, diamond,
kagome lattices were discussed in \cite{AIM14}.
Laplacians on periodic graphs with
non-compact perturbations and the stability of their
essential spectrum were considered in \cite{SS15}. The spectrum of Laplacians on the
lattice $\Z^d$ with pendant edges was studied in \cite{S13}. In the
paper \cite{KS16a}, the authors considered Schr\"odinger operators
with periodic potentials on periodic discrete graphs perturbed by
so-called guided potentials, which are periodic in some directions
and finitely supported in others. They described some properties
of the additional guided spectrum. We remark that the case  of
guided potentials is simpler than the case of periodic graphs with
guides and helps us to understand better the properties of the guided spectrum in the case of periodic graphs with guides. It is important that in the case of guided potentials all operators act in the same space. But in the case of periodic graphs with guides this is not true. Note that line
defects on the lattice were considered in \cite{C12}, \cite{Ku14},
\cite{Ku16}, \cite{OA12}.

Scattering theory for self-adjoint Schr\"odinger operators with
decreasing potentials  was investigated in \cite{BS99}, \cite{IK12}
(for the lattice) and in \cite{PR16} (for  periodic graphs). Inverse
scattering theory with finitely supported potentials was considered
in \cite{IK12} for the case of the lattice $\Z^d$ and in \cite{A12}
for the case of the hexagonal lattice. The absence of eigenvalues
embedded in the essential spectrum of the operators was discussed in
\cite{IM14}, \cite{V14}. Trace formulae and global eigenvalues
estimates for Schr\"odinger operators with complex decaying
potentials on the lattice were obtained in \cite{KL16}. The
Cwikel-Lieb-Rosenblum type bound for the discrete Schr{\"o}dinger
operator on $\Z^d$ was computed in \cite{Ka08}, \cite{RS09}.
Finally, we note that different properties of Schr\"odinger
operators on graphs  were considered in \cite{G15}, \cite{Sh98}.

\section {\lb{Sec1.1}Main results}
\setcounter{equation}{0}

\subsection{The unperturbed case: periodic graphs}
We define the infinite \emph{fundamental graph} $\cC_0$ of the
$\Z^{\wt d}$-periodic graph $\G_0$ by
$$
\cC_0=(V^c_0,\cE^c_0)=\G_0/\Z^d,\qqq V_0^c=V_0/\Z^d,\qqq
\cE_0^c=\cE_0/\Z^d, \qqq d<\wt d,
$$
where $V_0^c$ is its vertex set and $\cE_0^c$ is its set of
unoriented edges. Remark that the graph $\cC_0$ is a graph on the
cylinder $\R^{\wt d}/\Z^d$ and is $\Z^{\wt d-d}$-periodic. We also
call the fundamental graph $\cC_0$ \emph{a discrete cylinder} or just
\emph{a cylinder}. We identify the vertices of the cylinder $\cC_0$
with the vertices of the periodic graph $\G_0$ from the strip
$\cS=[0,1)^d\ts\R^{\wt d-d}$. We will call this infinite vertex set
\emph{a fundamental vertex set of $\G_0$} and denote it by the same
symbol $V^c_0$:
\[\lb{fvs}
V^c_0=V_0\cap\cS,\qqq \cS=[0,1)^d\ts\R^{\wt d-d}.
\]
Edges of the periodic graph $\G_0$ connecting the vertices from  the
fundamental vertex set $V^c_0$ with the vertices from $V_0\sm V^c_0$
will be called \emph{bridges}. Bridges always exist and provide the
connectivity of the periodic graph. The set of all bridges of the
graph $\G_0$ we denote by $\cB$.

\subsection{The perturbed case: periodic graphs with guides}
We define the \emph{union} of two graphs $G_0=(\cV_0,\cE_0)$ and
$G_1=(\cV_1,\cE_1)$ is a graph $G$ given by
$$
G=G_0\cup G_1=(\wt V,\wt \cE\,)\qqq \wt V=\cV_0\cup \cV_1,\qqq
\wt\cE=\cE_0\cup\cE_1.
$$

Now we define a periodic graph with guides. Let $\G_1=(V_1,\cE_1)$
be a finite graph, possibly not connected, such that all vertices of
$\G_1$ are contained in the strip $\cS$ and the graph $\G_0\cup\G_1$
is connected.
 A $\Z^d$-periodic graph
\[\lb{guid}
\G_1^g=\bigcup_{m\in\Z^d}(\G_1+m)
\]
will be called \emph{a guide} with the fundamental graph $\G_1$. We
define a perturbed graph $\G$ as a union of the unperturbed periodic
graph $\G_0$ and the perturbation $\G_1^g$:
\[
\lb{pgwg}
\G=\G_0\cup\G_1^g.
\]
We will call the graph $\G$ \emph{a periodic graph with a guide
$\G_1^g$} or \emph{a perturbed graph}.

Due to the definition \er{pgwg} of the perturbed graph
$\G$, the \emph{perturbed} cylinder $\cC=\G/\Z^d=(V^c,\cE^c)$ for
$\G$ is a union of the unperturbed cylinder $\cC_0=(V^c_0,\cE^c_0)$
and the finite graph $\G_1=(V_1,\cE_1)$:
\[
\cC=\cC_0\cup \G_1,
\]
i.e.,
\[\lb{uni0}
V^c=V^c_0\cup V_1,\qqq \cE^c=\cE^c_0\cup\cE_1, \qqq \cA^c=\cA^c_0\cup\cA_1,
\]
where $\cA^c$, $\cA^c_0$ and $\cA_1$ are the sets of all doubled
oriented  edges of $\cC$, $\cC_0$ and $\G_1$, respectively.

\setlength{\unitlength}{1.0mm}
\begin{figure}[h]
\centering

\unitlength 0.7mm % = 2.845pt
\linethickness{0.4pt}
\ifx\plotpoint\undefined\newsavebox{\plotpoint}\fi % GNUPLOT compatibility
\begin{picture}(250,70)(0,0)

\put(0,10){(\emph{a})}
\put(0,64){$\dL^2$}
\put(32,64){$\cS$}
%\put(77,38.5){$\Rightarrow$}

\put(0,20){\line(1,0){70.00}}
\put(0,30){\line(1,0){70.00}}
\put(0,40){\line(1,0){70.00}}
\put(0,50){\line(1,0){70.00}}
\put(0,60){\line(1,0){70.00}}

\bezier{60}(31,10)(31,40)(31,70)
\bezier{60}(32,10)(32,40)(32,70)
\bezier{60}(33,10)(33,40)(33,70)
\bezier{60}(34,10)(34,40)(34,70)
\bezier{60}(35,10)(35,40)(35,70)
\bezier{60}(36,10)(36,40)(36,70)
\bezier{60}(37,10)(37,40)(37,70)
\bezier{60}(38,10)(38,40)(38,70)
\bezier{60}(39,10)(39,40)(39,70)
\bezier{60}(41,10)(41,40)(41,70)
\bezier{60}(42,10)(42,40)(42,70)
\bezier{60}(43,10)(43,40)(43,70)
\bezier{60}(44,10)(44,40)(44,70)
\bezier{60}(45,10)(45,40)(45,70)
\bezier{60}(46,10)(46,40)(46,70)
\bezier{60}(47,10)(47,40)(47,70)
\bezier{60}(48,10)(48,40)(48,70)
\bezier{60}(49,10)(49,40)(49,70)
\bezier{60}(50,10)(50,40)(50,70)

\put(10,10){\line(0,1){60.00}}
\put(20,10){\line(0,1){60.00}}
\put(30,10){\line(0,1){60.00}}
\put(40,10){\line(0,1){60.00}}
\put(50,10){\line(0,1){60.00}}
\put(60,10){\line(0,1){60.00}}

\put(10,20){\circle{1}}
\put(20,20){\circle{1}}
\put(30,20){\circle*{2}}
\put(40,20){\circle*{2}}
\put(50,20){\circle{1}}
\put(60,20){\circle{1}}

\put(10,30){\circle{1}}
\put(20,30){\circle{1}}
\put(30,30){\circle*{2}}
\put(40,30){\circle*{2}}
\put(30,30){\vector(1,0){20.00}}
\put(30,30){\vector(0,1){10.00}}
\put(25,25){$O$}
\put(45,26){$a_1$}
\put(24,36.5){$a_2$}
\put(50,30){\circle{1}}
\put(60,30){\circle{1}}

\put(10,40){\circle{1}}
\put(20,40){\circle{1}}
\put(30,40){\circle*{2}}
\put(40,40){\circle*{2}}
\put(50,40){\circle{1}}
\put(60,40){\circle{1}}

\put(10,50){\circle{1}}
\put(20,50){\circle{1}}
\put(30,50){\circle*{2}}
\put(40,50){\circle*{2}}
\put(50,50){\circle{1}}
\put(60,50){\circle{1}}

\put(10,60){\circle{1}}
\put(20,60){\circle{1}}
\put(30,60){\circle*{2}}
\put(40,60){\circle*{2}}
\put(50,60){\circle{1}}
\put(60,60){\circle{1}}

%***************************
\put(80,10){\line(0,1){60.0}}
\put(90,10){\line(0,1){60.0}}
\bezier{60}(81,10)(81,40)(81,70)
\bezier{60}(82,10)(82,40)(82,70)
\bezier{60}(83,10)(83,40)(83,70)
\bezier{60}(84,10)(84,40)(84,70)
\bezier{60}(85,10)(85,40)(85,70)
\bezier{60}(86,10)(86,40)(86,70)
\bezier{60}(87,10)(87,40)(87,70)
\bezier{60}(88,10)(88,40)(88,70)
\bezier{60}(89,10)(89,40)(89,70)
\bezier{60}(91,10)(91,40)(91,70)
\bezier{60}(92,10)(92,40)(92,70)
\bezier{60}(93,10)(93,40)(93,70)
\bezier{60}(94,10)(94,40)(94,70)
\bezier{60}(95,10)(95,40)(95,70)
\bezier{60}(96,10)(96,40)(96,70)
\bezier{60}(97,10)(97,40)(97,70)
\bezier{60}(98,10)(98,40)(98,70)
\bezier{60}(99,10)(99,40)(99,70)
\bezier{60}(100,10)(100,40)(100,70)

%\multiput(100,11)(0,4){15}{\line(0,1){2}}

\put(80,20){\line(1,0){20.0}}
\put(80,30){\line(1,0){20.0}}
\put(80,40){\line(1,0){20.0}}
\put(80,50){\line(1,0){20.0}}
\put(80,60){\line(1,0){20.0}}

\put(80,20){\circle*{1.5}}
\put(90,20){\circle*{1.5}}
\put(100,20){\circle{1.5}}
\put(80,30){\circle*{1.5}}
\put(90,30){\circle*{1.5}}
\put(100,30){\circle{1.5}}
\put(80,40){\circle*{1.5}}
\put(90,40){\circle*{1.5}}
\put(100,40){\circle{1.5}}
\put(80,50){\circle*{1.5}}
\put(90,50){\circle*{1.5}}
\put(100,50){\circle{1.5}}
\put(80,60){\circle*{1.5}}
\put(90,60){\circle*{1.5}}
\put(100,60){\circle{1.5}}

\put(73,64){$\cC_0$}
\put(71,10){(\emph{b})}
%\put(118,38.5){$\Rightarrow$}
%*******************************

\put(220,10){\line(0,1){60.0}}
\put(230,10){\line(0,1){60.0}}
\bezier{60}(221,10)(221,40)(221,70)
\bezier{60}(222,10)(222,40)(222,70)
\bezier{60}(223,10)(223,40)(223,70)
\bezier{60}(224,10)(224,40)(224,70)
\bezier{60}(225,10)(225,40)(225,70)
\bezier{60}(226,10)(226,40)(226,70)
\bezier{60}(227,10)(227,40)(227,70)
\bezier{60}(228,10)(228,40)(228,70)
\bezier{60}(229,10)(229,40)(229,70)
\bezier{60}(231,10)(231,40)(231,70)
\bezier{60}(232,10)(232,40)(232,70)
\bezier{60}(233,10)(233,40)(233,70)
\bezier{60}(234,10)(234,40)(234,70)
\bezier{60}(235,10)(235,40)(235,70)
\bezier{60}(236,10)(236,40)(236,70)
\bezier{60}(237,10)(237,40)(237,70)
\bezier{60}(238,10)(238,40)(238,70)
\bezier{60}(239,10)(239,40)(239,70)
\bezier{60}(240,10)(240,40)(240,70)

%\multiput(240,11)(0,4){15}{\line(0,1){2}}

\put(220,20){\line(1,0){20.0}}
\put(220,30){\line(1,0){20.0}}
\put(220,40){\line(1,0){20.0}}
\put(220,50){\line(1,0){20.0}}
\put(220,60){\line(1,0){20.0}}

\put(220,20){\circle*{1.5}}
\put(230,20){\circle*{1.5}}
\put(240,20){\circle{1.5}}
\put(220,30){\circle*{1.5}}
\put(230,30){\circle*{1.5}}
\put(240,30){\circle{1.5}}
\put(220,40){\circle*{1.5}}
\put(230,40){\circle*{1.5}}
\put(240,40){\circle{1.5}}
\put(220,50){\circle*{1.5}}
\put(230,50){\circle*{1.5}}
\put(240,50){\circle{1.5}}
\put(220,60){\circle*{1.5}}
\put(230,60){\circle*{1.5}}
\put(240,60){\circle{1.5}}

\put(230,40){\line(1,2){3.00}}
\put(233,46){\circle*{1.5}}

\bezier{100}(220,30)(222.5,32.5)(225,35)
\bezier{100}(230,30)(227.5,32.5)(225,35)
\put(225,35){\circle*{1.5}}

\put(210,10){(\emph{e})}
\put(213,64){$\cC$}
%\put(208,38.5){$\Rightarrow$}

%**********************************

\put(140,20){\line(1,0){70.00}}
\put(140,30){\line(1,0){70.00}}
\put(140,40){\line(1,0){70.00}}
\put(140,50){\line(1,0){70.00}}
\put(140,60){\line(1,0){70.00}}

\put(150,10){\line(0,1){60.00}}
\put(160,10){\line(0,1){60.00}}
\put(170,10){\line(0,1){60.00}}
\put(180,10){\line(0,1){60.00}}
\put(190,10){\line(0,1){60.00}}
\put(200,10){\line(0,1){60.00}}

\bezier{60}(171,10)(171,40)(171,70)
\bezier{60}(172,10)(172,40)(172,70)
\bezier{60}(173,10)(173,40)(173,70)
\bezier{60}(174,10)(174,40)(174,70)
\bezier{60}(175,10)(175,40)(175,70)
\bezier{60}(176,10)(176,40)(176,70)
\bezier{60}(177,10)(177,40)(177,70)
\bezier{60}(178,10)(178,40)(178,70)
\bezier{60}(179,10)(179,40)(179,70)
\bezier{60}(181,10)(181,40)(181,70)
\bezier{60}(182,10)(182,40)(182,70)
\bezier{60}(183,10)(183,40)(183,70)
\bezier{60}(184,10)(184,40)(184,70)
\bezier{60}(185,10)(185,40)(185,70)
\bezier{60}(186,10)(186,40)(186,70)
\bezier{60}(187,10)(187,40)(187,70)
\bezier{60}(188,10)(188,40)(188,70)
\bezier{60}(189,10)(189,40)(189,70)

\put(160,40){\line(1,2){3.00}}
\put(163,46){\circle{1}}

\bezier{100}(150,30)(152.5,32.5)(155,35)
\bezier{100}(160,30)(157.5,32.5)(155,35)
\put(155,35){\circle{1}}

\bezier{100}(170,30)(172.5,32.5)(175,35)
\bezier{100}(180,30)(177.5,32.5)(175,35)
\put(175,35){\circle*{2}}

\bezier{100}(190,30)(192.5,32.5)(195,35)
\bezier{100}(200,30)(197.5,32.5)(195,35)
\put(195,35){\circle{1}}

\put(180,40){\line(1,2){3.00}}
\put(183,46){\circle*{2}}

\put(200,40){\line(1,2){3.00}}
\put(203,46){\circle{1}}

\put(150,20){\circle{1}}
\put(160,20){\circle{1}}
\put(170,20){\circle*{2}}
\put(180,20){\circle*{2}}
\put(190,20){\circle{1}}
\put(200,20){\circle{1}}

\put(150,30){\circle{1}}
\put(160,30){\circle{1}}
\put(170,30){\circle*{2}}
\put(180,30){\circle*{2}}
\put(190,30){\circle{1}}
\put(200,30){\circle{1}}
\put(170,30){\vector(1,0){20.00}}
\put(170,30){\vector(0,1){10.00}}
\put(165,25){$O$}
\put(185,26){$a_1$}
\put(164,36.5){$a_2$}

\put(150,40){\circle{1}}
\put(160,40){\circle{1}}
\put(170,40){\circle*{2}}
\put(180,40){\circle*{2}}
\put(190,40){\circle{1}}
\put(200,40){\circle{1}}

\put(150,50){\circle{1}}
\put(160,50){\circle{1}}
\put(170,50){\circle*{2}}
\put(180,50){\circle*{2}}
\put(190,50){\circle{1}}
\put(200,50){\circle{1}}

\put(150,60){\circle{1}}
\put(160,60){\circle{1}}
\put(170,60){\circle*{2}}
\put(180,60){\circle*{2}}
\put(190,60){\circle{1}}
\put(200,60){\circle{1}}

\put(140,64){$\G$}
\put(172,64){$\cS$}
\put(140,10){(\emph{d})}

%***************************
\put(115,64){$\G_1$}
\put(115,10){(\emph{c})}
\bezier{100}(110,30)(112.5,32.5)(115,35)
\bezier{100}(120,30)(117.5,32.5)(115,35)
\put(115,35){\circle*{1}}
\put(110,30){\circle*{1}}
\put(120,30){\circle*{1}}

\put(120,40){\line(1,2){3.00}}
\put(123,46){\circle*{1}}
\put(120,40){\circle*{1}}
\put(110,30){\vector(1,0){20.00}}
\put(110,30){\vector(0,1){10.00}}
\put(106,25){$O$}
\put(125,26){$a_1$}
\put(104,37.0){$a_2$}
\end{picture}

\vspace{-0.5cm} \caption{\footnotesize  \emph{a}) The square lattice $\dL^2$; the vertices from the set $V_0^c$ are big black points; the strip $\cS$ is shaded;\; \emph{b}) the unperturbed cylinder $\cC_0=\dL^2/\Z$ (the edges of the strip are identified); \; \emph{c}) a perturbation $\G_1$ with two connected components; \; \emph{d}) the perturbed square lattice $\G=\dL^2\cup\G_1^g$; \; \emph{e}) the perturbed cylinder $\cC=\G/\Z$ (the edges of the strip are identified).} \label{PGWG}
\end{figure}
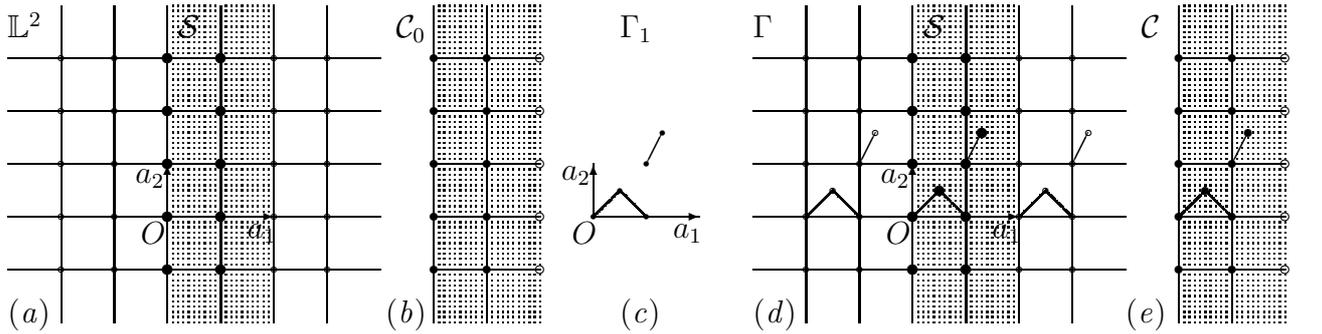

\textbf{Example.} For the square lattice $\dL^2$ with the periods
$a_1,a_2$,  see Fig.\ref{PGWG}.\emph{a}, the unperturbed cylinder
$\cC_0=\dL^2/\Z=(V_0^c,\cE_0^c)$ is shown in
Fig.\ref{PGWG}.\emph{b}. The vertices from the set $V_0^c$ are big
black points in Fig.\ref{PGWG}.\emph{a}. A perturbation $\G_1$, the
perturbed square lattice $\G=\dL^2\cup\G_1^g$ and the perturbed
cylinder $\cC=\G/\Z$ are shown in
Fig.\ref{PGWG}.\emph{c},\emph{d},\emph{e}.

\subsection{Floquet decomposition and the spectrum of perturbed Laplacians}
We describe the basic spectral  properties of Laplacians on periodic
graphs with guides.

\begin{proposition}\label{TSDI}
i) The Laplacian $\D$ on a perturbed graph $\G$ has the
following  decomposition into a constant fiber direct integral for
some unitary operator $U: \ell^2(V)\to \mH$:
\[
\lb{raz}
\begin{aligned}
& \mH=\int^\oplus_{\T^d}\ell^2(V^c){d\vt\/(2\pi)^d}\,,\qqq U\D
U^{-1}=\int^\oplus_{\T^d}\D(\vt){d\vt\/(2\pi)^d}\,,
\end{aligned}
\]
where $\T^d=\R^d/(2\pi\Z)^d$ and the fiber Laplacian $\D(\vt)$
 on the fiber space $\ell^2(V^c)$ is given by
\[
\label{l2.15''}
 \big(\D(\vt)f\big)(v)=\sum_{\be=(v,\,u)\in\cA^c}
 \big(f(v)-e^{i\lan\t(\be),\,\vt\ran}f(u)\big), \qqq
 v\in V^c,\qqq f\in\ell^2(V^c).
\]
Here $\t(\be)\in\Z^d$ is the index of the edge $\be\in\cA^c$ defined
by \er{in},\er{inf}, $V^c$ and $\cA^c$ are the vertex set and  the
set of oriented edges of the cylinder $\cC$, respectively;
$\lan\cdot\,,\cdot\ran$ is the inner product in~$\R^d$.

ii) For each $\vt\in\T^d$ the spectrum of the fiber operator
$\D(\vt)$ has the form
\[\lb{strs}
\s\big(\D(\vt)\big)=\s_{ac}\big(\D(\vt)\big)\cup\s_{fb}\big(\D(\vt)\big)
\cup\s_{p}\big(\D(\vt)\big),
\]
\[\lb{strs1}
\s_{ac}\big(\D(\vt)\big)=\s_{ac}\big(\D_0(\vt)\big),\qqq
\s_{fb}\big(\D(\vt)\big) =\s_{fb}\big(\D_0(\vt)\big),
\]
where $\D_0(\vt)$ is the fiber operator for the unperturbed
Laplacian $\D_0$  on the periodic graph $\G_0$,
$\s_{p}\big(\D(\vt)\big)$ is the set of all eigenvalues of $\D(\vt)$
of finite multiplicity given by
\[
\label{eq.3'}
\l_{N_\vt}(\vt)\le\ldots\le\l_2(\vt)\le\l_1(\vt),\qq
N_\vt\leq p:=\rank\D_1=\n_1-c_{\G_1},\qq \n_1=\#V_1,
\]
$\D_1$ is the Laplacian on the finite graph $\G_1=(V_1,\cE_1)$ and $c_{\G_1}$ is
 the number of connected components of $\G_1$, $\# A$ denotes
  the number of all elements of the set $A$.
\end{proposition}

\no \textbf{Remark.} The fiber Laplacian $\D(\vt)$, $\vt\in\T^d$,
can be considered as a magnetic Laplacian on the cylinder $\cC$ (see
\cite{HS99}, \cite{KS17}).

\medskip

Proposition \ref{TSDI} and standard arguments (see Theorem XIII.85
in \cite{RS78}) describe the spectrum of the Laplacian $\D$ on a
perturbed graph $\G$. Since $\D(\vt)$ is self-adjoint and
real analytic on the torus $\T^d=\R^d/(2\pi\Z)^d$, each
$\l_j(\cdot)$ is a real and piecewise analytic function on $\T^d$
and creates the \emph{guided band} $\gs_j(\D)$ given by
\[\lb{sgSo}
\gs_j(\D)=[\l_j^-,\l_j^+]=\l_j(\T^d), \qqq j=1,\ldots,N,\qqq
N=\max\limits_{\vt\in\T^d}N_\vt\leq p.
\]
Thus, the spectrum of the Laplacian $\D$ on the perturbed 
graph  $\G$ has the form
$$
\s(\D)=\bigcup_{\vt\in\T^d}\s\big(\D(\vt)\big)=\s(\D_0)\cup\gs(\D),
$$
where $\s(\D_0)$ is defined by \er{sH0} and
\[\lb{dgsp}
\gs(\D)=\bigcup\limits_{\vt\in\T^d}\s_{p}\big(\D(\vt)\big)
=\bigcup_{j=1}^N\gs_j(\D)=\gs_{ac}(\D)\cup\gs_{fb}(\D),
\]
$\gs_{ac}(\D)$ and $\gs_{fb}(\D)$ are the absolutely continuous part
and the flat band part of the guided spectrum $\gs(\D)$,
respectively. An open interval between two neighboring
non-degenerate bands is called \emph{a spectral gap}. The guided
spectrum $\gs(\D)$  may partly lie above the spectrum of the
unperturbed operator $\D_0$, on the spectrum of $\D_0$ and in the gaps
of $\D_0$.

We formulate simple sufficient conditions for the existence of
guided  flat bands of the perturbed Laplacian $\D$.

\begin{proposition}\label{Pdis}
Let $\G_0$ be a periodic graph. Assume that $\z$ is an eigenvalue of
the Laplacian $\D_1$ on a finite graph $\G_1$ with an
eigenfunction $f\in\ell^2(V_1)$ equal to zero on $V_0^c\cap V_1$,
i.e.,
\[\lb{Dico}
f(v)=0,\qqq \forall v\in V_{01}=V_0^c\cap V_1.
\]
Then $\{\z\}$ is a guided flat band of the Laplacian $\D$ on the
perturbed graph $\G=\G_0\cup\G_1^g$.
\end{proposition}

\setlength{\unitlength}{1.0mm}
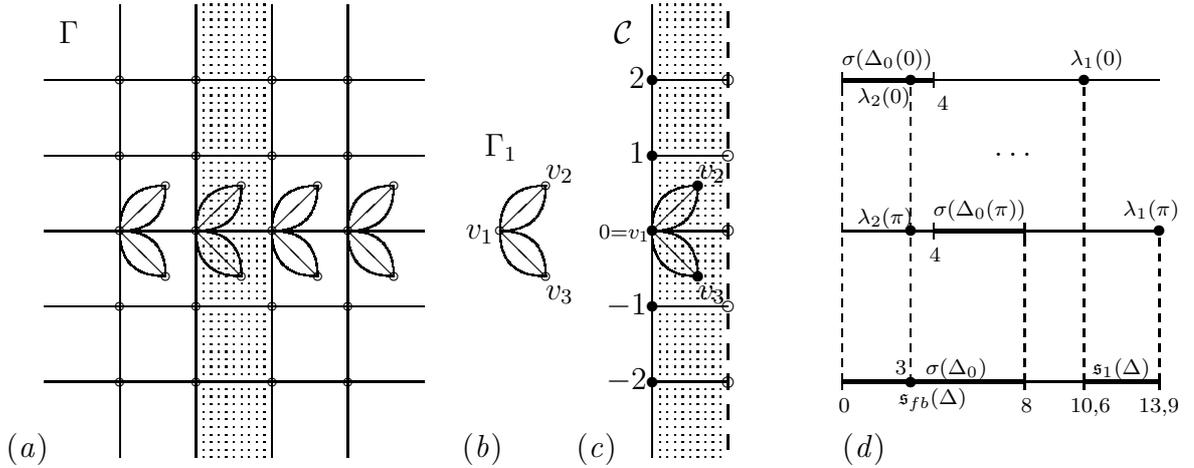
\begin{figure}[h]
\centering

\unitlength 1.0mm % = 2.845pt
\linethickness{0.4pt}
\ifx\plotpoint\undefined\newsavebox{\plotpoint}\fi % GNUPLOT compatibility
\begin{picture}(180,70)(0,0)
\put(5,10){(\emph{a})}
\put(12,65){$\G$}

\put(10,20){\line(1,0){50.00}}
\put(10,30){\line(1,0){50.00}}
\put(10,40){\line(1,0){50.00}}
\put(10,50){\line(1,0){50.00}}
\put(10,60){\line(1,0){50.00}}

\put(20,10){\line(0,1){60.00}}
\put(30,10){\line(0,1){60.00}}
\put(40,10){\line(0,1){60.00}}
\put(50,10){\line(0,1){60.00}}

\bezier{60}(31,10)(31,40)(31,70)
\bezier{60}(32,10)(32,40)(32,70)
\bezier{60}(33,10)(33,40)(33,70)
\bezier{60}(34,10)(34,40)(34,70)
\bezier{60}(35,10)(35,40)(35,70)
\bezier{60}(36,10)(36,40)(36,70)
\bezier{60}(37,10)(37,40)(37,70)
\bezier{60}(38,10)(38,40)(38,70)
\bezier{60}(39,10)(39,40)(39,70)

%**********************************
\put(20,40){\line(1,1){6.00}}
\put(26,46){\circle{1}}
\bezier{100}(20,40)(20,46)(26,46)
\bezier{100}(20,40)(26,40)(26,46)

\put(30,40){\line(1,1){6.00}}
\put(36,46){\circle{1}}
\bezier{100}(30,40)(30,46)(36,46)
\bezier{100}(30,40)(36,40)(36,46)

\put(40,40){\line(1,1){6.00}}
\put(46,46){\circle{1}}
\bezier{100}(40,40)(40,46)(46,46)
\bezier{100}(40,40)(46,40)(46,46)

\put(50,40){\line(1,1){6.00}}
\put(56,46){\circle{1}}
\bezier{100}(50,40)(50,46)(56,46)
\bezier{100}(50,40)(56,40)(56,46)

%***************************
\put(20,40){\line(1,-1){6.00}}
\put(26,34){\circle{1}}
\bezier{100}(20,40)(20,34)(26,34)
\bezier{100}(20,40)(26,40)(26,34)

\put(30,40){\line(1,-1){6.00}}
\put(36,34){\circle{1}}
\bezier{100}(30,40)(30,34)(36,34)
\bezier{100}(30,40)(36,40)(36,34)

\put(40,40){\line(1,-1){6.00}}
\put(46,34){\circle{1}}
\bezier{100}(40,40)(40,34)(46,34)
\bezier{100}(40,40)(46,40)(46,34)

\put(50,40){\line(1,-1){6.00}}
\put(56,34){\circle{1}}
\bezier{100}(50,40)(50,34)(56,34)
\bezier{100}(50,40)(56,40)(56,34)

%****************************
\put(20,20){\circle{1}}
\put(30,20){\circle{1}}
\put(40,20){\circle{1}}
\put(50,20){\circle{1}}

\put(20,30){\circle{1}}
\put(30,30){\circle{1}}
\put(40,30){\circle{1}}
\put(50,30){\circle{1}}

\put(20,40){\circle{1}}
\put(30,40){\circle{1}}
\put(40,40){\circle{1}}
\put(50,40){\circle{1}}

\put(20,50){\circle{1}}
\put(30,50){\circle{1}}
\put(40,50){\circle{1}}
\put(50,50){\circle{1}}

\put(20,60){\circle{1}}
\put(30,60){\circle{1}}
\put(40,60){\circle{1}}
\put(50,60){\circle{1}}

%****************************
\put(65,10){(\emph{b})}
\put(68,50){$\G_1$}
\put(70,40){\line(1,1){6.00}}
\put(76,46){\circle{1}}
\bezier{100}(70,40)(70,46)(76,46)
\bezier{100}(70,40)(76,40)(76,46)
\put(70,40){\circle{1}}
\put(70,40){\line(1,-1){6.00}}
\put(76,34){\circle{1}}
\bezier{100}(70,40)(70,34)(76,34)
\bezier{100}(70,40)(76,40)(76,34)
\put(65.5,39.0){$v_1$}
\put(76,31.0){$v_3$}
\put(76,47.0){$v_2$}
%***************************
\put(85,65){$\cC$}
\put(90,40){\line(1,1){6.00}}
\put(96,46){\circle*{1.5}}
\bezier{100}(90,40)(90,46)(96,46)
\bezier{100}(90,40)(96,40)(96,46)

\put(90,40){\line(1,-1){6.00}}
\put(96,34){\circle*{1.5}}
\bezier{100}(90,40)(90,34)(96,34)
\bezier{100}(90,40)(96,40)(96,34)

\put(90,10){\line(0,1){60.0}}
\multiput(100,11)(0,4){15}{\line(0,1){2}}

\bezier{60}(91,10)(91,40)(91,70)
\bezier{60}(92,10)(92,40)(92,70)
\bezier{60}(93,10)(93,40)(93,70)
\bezier{60}(94,10)(94,40)(94,70)
\bezier{60}(95,10)(95,40)(95,70)
\bezier{60}(96,10)(96,40)(96,70)
\bezier{60}(97,10)(97,40)(97,70)
\bezier{60}(98,10)(98,40)(98,70)
\bezier{60}(99,10)(99,40)(99,70)

\put(90,20){\line(1,0){10.0}}
\put(90,30){\line(1,0){10.0}}
\put(90,40){\line(1,0){10.0}}
\put(90,50){\line(1,0){10.0}}
\put(90,60){\line(1,0){10.0}}

\put(90,20){\circle*{1.5}}
\put(100,20){\circle{1.5}}

\put(90,30){\circle*{1.5}}
\put(100,30){\circle{1.5}}

\put(90,40){\circle*{1.5}}
\put(100,40){\circle{1.5}}

\put(90,50){\circle*{1.5}}
\put(100,50){\circle{1.5}}

\put(90,60){\circle*{1.5}}
\put(100,60){\circle{1.5}}

\put(80,10){(\emph{c})}
\put(84,19.0){$-2$}
\put(84,29.0){$-1$}
\put(83,39.0){$\scriptstyle 0=v_1$}
\put(87,49.0){$1$}
\put(87,59.0){$2$}

\put(96,31.0){$v_3$}

\put(96,47.0){$v_2$}

%***************************
\put(115,10){(\emph{d})}
\multiput(115,20)(0,2){20}{\line(0,1){1}}
\put(115,20){\line(1,0){41.7}}
\put(115,19){\line(0,1){2.00}}
\put(115,20.2){\line(1,0){24.00}}
\put(115,19.8){\line(1,0){24.00}}
\put(139,19){\line(0,1){2.00}}
\multiput(139,20)(0,2){10}{\line(0,1){1}}
\put(114.5,16.0){$\scriptstyle 0$}
\put(138.5,16.0){$\scriptstyle 8$}
\put(124,20){\circle*{1.5}}
\multiput(124,20)(0,2){20}{\line(0,1){1}}
\put(122.5,17){$\scriptstyle \gs_{fb}(\D)$}
\multiput(146.8,20)(0,2){20}{\line(0,1){1}}
\put(146.8,19){\line(0,1){2.00}}
\put(146.8,19.8){\line(1,0){9.9}}
\put(146.8,20.2){\line(1,0){9.9}}
\put(156.7,19){\line(0,1){2.00}}
\multiput(156.7,20)(0,2){10}{\line(0,1){1}}
\put(145,16){$\scriptstyle 10{,}6$}
\put(154,16){$\scriptstyle 13{,}9$}
\put(126,21){$\scriptstyle \s(\D_0)$}
\put(122,21){$\scriptstyle 3$}
\put(148,21){$\scriptstyle \gs_1(\D)$}

\put(115,60){\line(1,0){41.7}}
\put(115,59){\line(0,1){2.00}}
\put(115,60.2){\line(1,0){12.00}}
\put(115,59.8){\line(1,0){12.00}}
\put(127,59){\line(0,1){2.00}}
%\put(114.5,56.0){$\scriptstyle 0$}
\put(127.5,56.0){$\scriptstyle 4$}
\put(124,60){\circle*{1.5}}
\put(117.0,57){$\scriptstyle \l_2(0)$}
\put(146.8,60){\circle*{1.5}}
%\put(145,56){$\scriptstyle 10{,}6$}
\put(115,62){$\scriptstyle \s(\D_0(0))$}
\put(145,62){$\scriptstyle \l_1(0)$}

\put(115,40){\line(1,0){41.7}}
\put(127,39){\line(0,1){2.00}}
\put(127,40.2){\line(1,0){12.00}}
\put(127,39.8){\line(1,0){12.00}}
\put(139,39){\line(0,1){2.00}}
\put(126.5,36.0){$\scriptstyle 4$}
%\put(138.0,36.0){$\scriptstyle 8$}
\put(124,40){\circle*{1.5}}
\put(117.0,41){$\scriptstyle \l_2(\pi)$}
\put(156.7,40){\circle*{1.5}}
%\put(154,36){$\scriptstyle 13{,}9$}
\put(127,42){$\scriptstyle \s(\D_0(\pi))$}
\put(152,42){$\scriptstyle \l_1(\pi)$}

\put(135,50){$\ldots$}
\end{picture}

\vspace{-0.5cm} \caption{\footnotesize  \emph{a}) The perturbed square lattice $\G=\dL^2\cup\G_1^g$; \; \emph{b}) The finite graph $\G_1$; \; \emph{c})~the~perturbed cylinder $\cC=\G/\Z$; \; \emph{d}) the spectra of the fiber Laplacians $\D(\vt)$ as $\vt=0,\pi$ and the Laplacian $\D$.} \label{SL2}
\end{figure}

\textbf{Example.} We consider the perturbed square lattice
$\G=\dL^2\cup\G_1^g$ shown in Fig.\ref{SL2}.\emph{a}. For each
$\vt\in\T=(-\pi,\pi]$ the spectrum of the fiber Laplacian $\D(\vt)$
has the form
$$
\s\big(\D(\vt)\big)=\s_{ac}\big(\D(\vt)\big)\cup\s_{p}\big(\D(\vt)\big), \qq
\s_{ac}\big(\D(\vt)\big)=\s\big(\D_0(\vt)\big)=[2-2\cos\vt,6-2\cos\vt],
$$
$\s_{p}\big(\D(\vt)\big)$ consists of two eigenvalues $\l_1(\vt)$
and $\l_2(\vt)=3$,  see Fig.\ref{SL2}.\emph{d}.

The spectrum of the Laplacian $\D$ on $\G$ has the form
$$
\s(\D)=\s(\D_0)\cup\gs(\D), \qqq \s(\D_0)=[0,8],
$$
where the guided spectrum $\gs(\D)$ is given by (see
Fig.\ref{SL2}.\emph{d} and  details in Proposition \ref{Prg3})
$$
\begin{aligned}
\gs(\D)=\gs_{ac}(\D)\cup\gs_{fb}(\D),\qqq \gs_{ac}(\D)=\gs_1(\D)
=\l_1(\T)\approx[10{,}6;13{,}9],\\
\gs_{fb}(\D)=\gs_2(\D)=\l_2(\T)=\{3\}.
\end{aligned}
$$
The Laplacian $\D_1$ on the finite graph $\G_1$ has the
eigenvalue  $\l=3$ with an eigenfunction $f$ such that $f(0)=0$. Then,
due to Proposition \ref{Pdis}, $\{3\}$ is a guided flat band of $\D$
on the perturbed square lattice $\G=\dL^2\cup\G_1^g$. Here we remark
that we do not know an example of a Schr\"odinger operator with a
guided potential on a periodic graph having a guided flat band.

\subsection{Estimates of guided bands.} We consider the guided bands from \er{sgSo}
 (or their parts) above the spectrum of the unperturbed Laplacian $\D_0$:
\[\lb{gs+}
\gs_j^o(\D)=\gs_j(\D)\cap[\vr,+\iy)\neq\varnothing,\qqq j=1,\ldots,N_g,\qqq N_g\leq N,
\]
recall that $\vr=\sup\s(\D_0)$. The Laplacian $\D_1$ on the finite
graph  $\G_1=(V_1,\cE_1)$ has the eigenvalue $0$ of multiplicity
$c_{\G_1}$ and $p$ positive eigenvalues $\z_j$ labeled by
\[\lb{ppp}
0<\z_p\le\ldots\le\z_2\le\z_1,\qqq p=\n_1-c_{\G_1},\qqq \n_1=\#V_1,
\]
 counting multiplicity,
where $c_{\G_1}$ is the number of connected components of the graph
$\G_1$.

Proposition \ref{TSDI} and the standard perturbation theory give the
estimates  of the position of the bands $\gs_j^o(\D)$ and their
number $N_g$ by (for more details see Corollary \ref{Tloc})
\[\lb{esbp0}
\gs_j^o(\D)\ss[\z_j,\z_j+\vr],\qqq N_g\geq\#\{j\in\N_p\,:\, \z_j>\vr\},
\]
 where $\#A$ is the number of elements of the set $A$. In
particular, this yields that if the eigenvalues of $\D_1$ satisfy
$\z_p>\vr$ and $\z_j-\z_{j+1}>\vr$ for all $j\in\N_{p-1}$, then the
guided spectrum of the Laplacian $\D$ consists of exactly $p$ guided
bands separated by gaps.

\medskip

In order to formulate our main result we define the set
$\cB^c=\cB/\Z^d$ of all bridges of the cylinder $\cC=(V^c,\cE^c)$
and the modified cylinder $\cC^\gm=(V^c,\cE^c\sm\cB^c)$, which is
obtained from $\cC$ by deleting all its bridges.  We consider the
Laplacian $\D^\gm$ defined by \er{DOLN} on the modified cylinder
$\cC^\gm$. This Laplacian $\D^\gm$ has at most $p=\rank\D_1$
eigenvalues $\wt\m_1\geq\wt\m_2\geq\ldots$\,. Define $\m_j$ by
\[\lb{evdm}
\m_j=\max\{\wt\m_j, \sup\s_{ess}(\D^\gm)\}, \qqq j=1,2,\ldots,p.
\]
We estimate the position of the guided bands $\gs_j^o(\D)$
defined by \er{gs+} in terms of the eigenvalues of the operator
$\D^\gm$ and the number of bridges on the cylinder~$\cC$.

\begin{theorem}\lb{Est}
i) Let $\D$ be the Laplacian on a perturbed graph $\G$ and let $\m_j$ be defined by \er{evdm}. Then
each guided  band $\gs_j^o(\D)$, $j=1,\ldots,N_g$, defined by
\er{gs+} satisfies
\[\lb{esbp1}
\gs_j^o(\D)\ss[\m_j,\m_j+2\b_+], \qqq \b_+=\max_{v\in V^c}\b_v,
\]
where  $\b_v$ is the number
of bridges  on $\cC$ starting at the vertex $v\in V^c$.

ii) Moreover,  for any $\ve>0$ there exists a perturbed graph $\G$  such that each non-degenerate guided band length
$|\gs_j^o(\D)|>2\b_+-\ve$, $j=1,\ldots,N_g$.
\end{theorem}

\no \textbf{Remark.} For most of graphs the number $\b_+=1$, then the guided band length $|\gs_j^o(\D)|\leq2$
for all $j=1,\ldots,N_g$, but for specific graphs $\b_+$ may be any given
integer number.

%2) Note that see Propositions
%\ref{Prg1}, \ref{Prg3}.

\medskip

Let $\G_t=(V_1,\cE_t)$ be a finite graph obtained from the graph
$\G_1=(V_1,\cE_1)$ considering each edge of $\G_1$ to have the
multiplicity $t\in\N$. We consider the Laplacian $\D_t$ acting on a
perturbed graph $\G=\G_0\cup\G_t^g$ and discuss the guided
spectrum of $\D_t$ for large $t$.

\begin{theorem}
\lb{T17} Let $\D_t$ be the Laplacian on the perturbed graph
$\G=\G_0\cup\G_t^g$, where $\G_0$ is any periodic graph and $t\in\N$
is large enough. Then the guided spectrum of the Laplacian $\D_t$
consists of exactly $p$ guided bands separated by gaps, where $p$ is
defined in \er{ppp},  and the following statements hold true:

i) Let $\z_j$ for some $j\in\N_p$ be a simple positive eigenvalue of
the Laplacian $\D_1$ on the graph $\G_1$ with a normalized
eigenfunction $f_j\in\ell^2(V_1)$.

$\bullet$ If $f_j=0$ on $V_{01}=V_0^c\cap V_1$, then $\{t\z_j\}$ is a guided flat band of the Laplacian $\D_t$.

$\bullet$ If $f_j\neq0$ on $V_{01}$, then the guided band
$\gs_j(\D_t)=[\l_j^-(t),\l_j^+(t)]$ satisfies
\[
\lb{Qt}
\begin{aligned}
&\l_j^\pm(t)=t\z_j+W_j^\pm+O(1/t),\\
&|\gs_j(\D_t)|=W_j^\bu+O(1/t),
\end{aligned}
\]
as $t\to \iy$, where
\[\lb{Dpm}
\begin{aligned}
&W_j^-=\min_{\vt\in\T^d}W_j(\vt),\qqq
W_j^+=\max_{\vt\in\T^d}W_j(\vt),\\
& W_j^\bu=W_j^+-W_j^-,\qqq W_j^\bu\leq2\b_{01},
\end{aligned}
\]
for some function $W_j$ defined by the formula \er{cjvt}. Here
$\b_{01}$  is the number of all oriented bridges connecting the
vertices from $V_{01}$ on the cylinder $\cC$.

ii) In particular, if the set $V_{01}$ consists of one vertex $v$, then
\[\lb{cov}
f^2_j(v)\b_{01}\leq W_j^\bu\leq 2f^2_j(v)\b_{01}.
\]
Moreover, $W_j^\bu=0$ iff $\b_{01}$=0.

iii) Let all positive eigenvalues $\z_j$, $j\in\N_p$, of $\D_1$ be distinct.
Then the Lebesgue measure $|\gs(\D_t)|$ of
the guided spectrum of the Laplacian $\D_t$ satisfies
\[
\lb{Qt1} |\gs(\D_t)|=\sum\limits_{j=1}^p W^\bu_j+O(1/t).
\]

iv) In particular, if there is no bridge connecting the vertices
from  the set $V_{01}$ on the cylinder $\cC$, then $W^\bu_j=0$ for
each $j\in\N_p$ and the second identity in \er{Qt} and the formula
\er{Qt1} take the form $|\gs_j(\D_t)|=O(1/t)$, $j\in\N_p$, and
$|\gs(\D_t)|=O(1/t)$, respectively.
\end{theorem}

\smallskip

Now we describe geometric properties of the guided spectrum for periodic graphs with specific guides.

\begin{corollary}\lb{TEg} Let $\G_0$ be a periodic graph with
an unperturbed cylinder $\cC_0=(V_0^c,\cE_0^c)$. Then the following
statements hold true.

i) For any constant $C>0$ there exists a finite graph $\G_t$, $t\in\N$,
such that the Lebesgue measure of the guided spectrum $\gs(\D)$ of
the perturbed Laplacian $\D$ on $\G=\G_0\cup\G_t^g$ satisfies
$|\gs(\D)|>C$ and all guided bands are non-degenerate.

ii) Let, in addition, there exist a vertex $v\in V_0^c$ such that
there is no bridge on $\cC_0$ starting at $v$. Then for any small
$\ve>0$ there exists a finite graph $\G_t$ such that the Lebesgue
measure of the guided spectrum $\gs(\D)$ of the perturbed Laplacian
$\D$ on $\G=\G_0\cup\G_t^g$ satisfies $|\gs(\D)|<\ve$.

iii) For any constant $\l_0>0$ there exists a finite graph $\G_t$
such  that the guided spectrum $\gs(\D)$ of the perturbed Laplacian
$\D$ on $\G=\G_0\cup\G_t^g$ satisfies
$\gs(\D)\cap(\l_0,+\iy)\neq\varnothing$.

iv) There exists a finite graph $\G_1$ such that the guided spectrum
$\gs(\D)$ of the perturbed Laplacian $\D$ on $\G=\G_0\cup\G_1^g$ has
a degenerate guided band.
\end{corollary}

Thus, roughly speaking, the guided spectrum can be any set above the
unperturbed  spectrum. Its Lebesgue measure can be arbitrarily large
or arbitrarily small.

We present the plan of our paper. In Section \ref{Sec2}  we
introduce the notion of edge indices and prove Proposition
\ref{TSDI} about the decomposition  of the Laplacian on periodic
graphs with guides into a constant fiber direct integral. In Section
\ref{Sec4} we prove Theorems \ref{Est}, \ref{T17} and Corollary
\ref{TEg}. Section \ref{Sec5} is devoted to properties of the guided
spectrum for  the square lattice with specific guides.

\section{Direct integral for Laplacians on periodic graphs with guides}
\setcounter{equation}{0} \lb{Sec2}
\subsection{Edge indices.}

In order to give a decomposition of Laplacians on periodic graphs
with guides   into a constant fiber direct integral with a precise
representation of fiber operators we need to define {\it an edge
index}. Recall that an edge index was introduced in \cite{KS14} and
it was important to study the spectrum of Laplacians and
Schr\"odinger operators on periodic graphs, since fiber operators
are expressed in terms of edge indices (see \er{l2.15''}).

For any
$v\in V$ the following unique representation holds true:
\[
\lb{Dv} v=v_0+[v], \qqq v_0\in V^c,\qqq [v]\in\Z^d,
\]
where $V^c$ is the fundamental vertex set of the graph
$\G=(V,\cE)$ defined by
\[\lb{fvsG}
V^c=V\cap\cS,\qqq \cS=[0,1)^d\ts\R^{\wt d-d}.
\]
In other words, each vertex $v$ can be obtained
from a vertex $v_0\in V^c$ by the shift by a vector $[v]\in\Z^d$.
For any
oriented edge $\be=(u,v)\in\cA$ we define {\bf the edge "index"}
$\t(\be)$ as the integer vector given by
\[
\lb{in}
\t(\be)=[v]-[u]\in\Z^d,
\]
where, due to \er{Dv}, we have
$$
u=u_0+[u],\qquad v=v_0+[v], \qquad u_0,v_0\in V^c,\qquad [u],[v]\in\Z^d.
$$
We note that edges connecting vertices from the fundamental vertex set $V^c$
have zero indices.

We define a surjection $\gf_\cA:\cA\rightarrow\cA^c=\cA/\Z^d$, which
map each edge to its equivalence class. If $\be$ is an oriented edge
of the graph $\G$, then there is an oriented edge
$\be_*=\gf_{\cA}(\be)$ on the cylinder $\cC=\G/\Z^d$. For the edge
$\be_*\in\cA^c$ we define the edge index $\t(\be_*)$ by
\[
\lb{inf}
\t(\be_*)=\t(\be).
\]
In other words, edge indices of the cylinder $\cC$  are induced by
edge indices of the graph $\G$. Edges with nonzero indices are
called \emph{bridges}. Edge indices, generally speaking, depend on
the choice of the coordinate origin $O$ and the periods
$a_1,\ldots,a_{\wt d}$ of the graph $\G$. But in a fixed  coordinate
system indices of the cylinder edges are uniquely determined by
\er{inf}, since
$$
\t(\be+m)=\t(\be),\qqq \forall\, (\be,m)\in\cA \ts \Z^d.
$$
We note that, due to the definition of periodic graphs with guides,
\[\lb{zin}
\t(\be)=0, \qqq \forall\,\be\in\cA_1.
\]

\subsection{Direct integrals.} We prove Proposition \ref{TSDI}  about
 the decomposition of Laplacians on periodic graphs with guides into a constant fiber direct integral.

\

\no \textbf{Proof of Proposition \ref{TSDI}.i)} Repeating  the
arguments from the proof of Theorem 1.1 in \cite{KS14} we obtain
\er{raz}, \er{l2.15''}, where the unitary operator
$U:\ell^2(V)\to\mH$ has the form
\[
\lb{5001} (Uf)(\vt,v)=\sum\limits_{m\in\Z^d}e^{-i\lan m,\vt\ran }
f(v+m), \qqq (\vt,v)\in \T^d\ts V^c, \qqq f\in \ell^2(V).
\]
The Hilbert space $\mH$ defined in \er{raz} is
equipped with the norm
$\|g\|^2_{\mH}=\int_{\T^d}\|g(\vt,\cdot)\|_{\ell^2(V^c)}^2\frac{d\vt
}{(2\pi)^d}$\,, where the function $g(\vt,\cdot)\in\ell^2(V^c)$ for almost all
$\vt\in\T^d$. \qq $\BBox$

In order to prove the next item of Proposition \ref{TSDI} we need the following lemma.

\begin{lemma}\lb{lcom}
Let $P$ and $P_1$ be the orthogonal projections of $\ell^2(V^c)$
onto the subspaces $\ell^2(V_0^c)$ and $\ell^2(V_1)$, respectively.
Then each fiber Laplacian $\D(\vt)$, $\vt\in\T^d$, defined by
\er{l2.15''} has the following decomposition:
\[\lb{dec1}
\D(\vt)=P\D_0(\vt)P+P_1\D_1P_1,
\]
where $\D_0(\vt)$ is the fiber operator for the unperturbed
Laplacian $\D_0$  on the periodic graph $\G_0$, $\D_1$ is the
Laplacian on the finite graph $\G_1=(V_1,\cE_1)$.
\end{lemma}

\no\textbf{Proof.} The Laplacian $\D_0$ on the unperturbed periodic
graph  $\G_0=(V_0,\cE_0)$ has a decomposition into a constant fiber
direct integral for some unitary operator $U_0: \ell^2(V_0)\to
\mH_0$:
\[
\lb{raz0}
\begin{aligned}
& \mH_0=\int^\oplus_{\T^d}\ell^2(V_0^c){d\vt\/(2\pi)^d}\,,\qqq U_0\D_0
U_0^{-1}=\int^\oplus_{\T^d}\D_0(\vt){d\vt\/(2\pi)^d}\,,
\end{aligned}
\]
where the fiber Laplacian $\D_0(\vt)$ acts on the fiber space
$\ell^2(V_0^c)$ and is given by
\[
\label{l2.15'}
 \big(\D_0(\vt)f_0\big)(v)=\sum_{\be=(v,\,u)\in\cA_0^c} \big(f_0(v)-e^{i\lan\t(\be),\,\vt\ran}f_0(u)\big), \qqq
 v\in V_0^c, \qqq f_0\in\ell^2(V_0^c).
\]
For each $f\in\ell^2(V^c)$ we have
$$
\begin{aligned}
&\big\lan(P\D_0(\vt)P+P_1\D_1P_1)f, f\big\ran_{V^c}= \lan
\D_0(\vt)Pf,Pf\ran_{V_0^c}+\lan \D_1P_1f,P_1f\ran_{V_1}\\&=
\sum_{v\in V_0^c}\big(\D_0(\vt)Pf\big)(v)\,\ol f(v)+
\sum_{v\in V_1}\big(\D_1P_1f\big)(v)\,\ol f(v)=\sum_{v\in V_0^c\sm V_{01}}\big(\D_0(\vt)Pf\big)(v)\,\ol f(v)\\
&+
\sum_{v\in V_1\sm V_{01}}\big(\D_1P_1f\big)(v)\,\ol f(v)
+\sum_{v\in V_{01}}\big((\D_0(\vt)P+\D_1P_1)f\big)(v)\,\ol f(v),
\end{aligned}
$$
where $\lan\cdot,\cdot\ran_{V}$ denotes the inner product in
$\ell^2(V)$. Substituting the definitions  \er{DOLN} and \er{l2.15'}
of the Laplacian $\D_1$ and the fiber Laplacian $\D_0(\vt)$ into
this formula and using the identities \er{uni0}, \er{zin} and \er{l2.15''}, we
obtain
$$
\begin{aligned}
&\big\lan(P\D_0(\vt)P+P_1\D_1P_1)f, f\big\ran_{V^c}=\\
&=\sum_{v\in V_0^c\sm V_{01}}\sum_{\be=(v,\,u)\in\cA_0^c}
\big(f(v)-e^{i\lan\t(\be),\,\vt\ran}f(u)\big)\,\ol f(v)
+\sum_{v\in V_1\sm V_{01}}\sum_{(v,u)\in\cA_1}\big(f(v)-f(u)\big)\,\ol f(v)\\
&+\sum_{v\in V_{01}}\bigg(\sum_{\be=(v,\,u)\in\cA_0^c} \big(f(v)-e^{i\lan\t(\be),\,
\vt\ran}f(u)\big)+\sum_{(v,u)\in \cA_1}\big(f(v)-f(u)\big)\bigg)\,\ol f(v)\\
&=\sum_{v\in V^c}\sum_{\be=(v,\,u)\in\cA^c} \big(f(v)-e^{i\lan\t(\be),\,\vt\ran}f(u)\big)\,\ol f(v)
=\sum_{v\in V^c}\big(\D(\vt)f\big)(v)\,\ol f(v)=\lan\D(\vt)f,f\ran_{V^c},
\end{aligned}
$$
which implies \er{dec1}. \qq $\BBox$

\

\no \textbf{Proof of Proposition \ref{TSDI}.ii)} For each $\vt\in\T^d$ the unperturbed
 fiber Laplacian $\D_0(\vt)$ is $\Z^{\wt d-d}$-periodic. Then, using standard arguments
  (see Theorem XIII.85 in \cite{RS78}), we obtain that the spectrum of $\D_0(\vt)$ has the form
$$
\s\big(\D_0(\vt)\big)=\s_{ac}\big(\D_0(\vt)\big)\cup
\s_{fb}\big(\D_0(\vt)\big).
$$
Since the graph $\G_1$ is finite, for each
$\l\in\s_{fb}\big(\D_0(\vt)\big)$ there exists a corresponding
eigenfunction $f\in\ell^2(V_0^c)$ with a finite support (see, e.g.,
Theorem 4.5.2 in \cite{BK13}) not intersecting with $V_1$. Due to
\er{dec1}, $(f,0)\in\ell^2(V^c)$ is  an eigenfunction of $\D(\vt)$
with the same finite support and the same eigenvalue $\l$.  Thus,
$\l\in\s_{fb}\big(\D(\vt)\big)$ and vice versa. Since the operator
$\D_1$ has finite rank $p$, where $p$ is defined in
\er{eq.3'}, for each $\vt\in\T^d$ the spectrum $\s\big(\D(\vt)\big)$
of the fiber Laplacian $\D(\vt)$ is given by \er{strs}, where
$\s_{ac}\big(\D(\vt)\big)$, $\s_{fb}\big(\D(\vt)\big)$ satisfy
\er{strs1} and $\s_{p}\big(\D(\vt)\big)$ consists of $N_\vt\leq p$
eigenvalues \er{eq.3'}. \qq $\BBox$

\section{Proof of the main results}
\setcounter{equation}{0}
\lb{Sec4}

In this section we prove Theorem \ref{Est} about the position of
guided  bands and Theorem \ref{T17} about the asymptotics of the
guided bands for guides with large multiplicity of their edges.  We
prove Corollary \ref{TEg} about geometric properties of the guided spectrum for
periodic graphs with specific guides.

\subsection{Estimates for the guided spectrum.} Denote by $m_\pm(\vt)$
 the upper and lower endpoints of the spectrum of the unperturbed fiber Laplacian $\D_0(\vt)$:
\[\lb{mpmm}
m_-(\vt)=
\inf\s\big(\D_0(\vt)\big),\qqq m_+(\vt)= \sup\s\big(\D_0(\vt)\big).
\]
Then \er{mm} yields
\[
\lb{mvt}
\min_{\vt\in\T^d}m_-(\vt)=0,\qqq \max_{\vt\in\T^d}m_+(\vt)=\vr.
\]

We need a simple estimate for eigenvalues of bounded self-adjoint operators \cite{RS78}: \emph{Let $A,B$ be bounded self-adjoint 
operators in a Hilbert space $\cH$ and let
$\l_j(A)=\max\big\{\wt\l_j(A),\sup\s_{ess}(A)\big\}$,
$j=1,2,\ldots$\,, where $\wt \l_1(A)\geq\wt \l_2(A)\geq\ldots$ are
the eigenvalues of $A$. Then}
\[\lb{wnf1}
\l_j(A)+\inf\s(B)\leq\l_j(A+B)\leq\l_j(A)+\sup\s(B), \qqq j=1,2,3,\ldots\,.
\]

The following simple corollary about the position of the guided
bands $\gs_j^o(\D)$  defined by \er{gs+} is a direct consequence of
Proposition \ref{TSDI}.

\begin{corollary}\lb{Tloc}
Let $\D$ be the Laplacian on a perturbed graph $\G$ and let
$\vr=\sup \s(\D_0)$. Then each guided band $\gs_j^o(\D)$,
$j=1,\ldots,N_g$, and their number $N_g$ satisfy
\[\lb{esbp}
\gs_j^o(\D)\ss[\z_j,\z_j+\vr],
\]
\[\lb{esNg}
N_g\geq\#\{j\in\N_p\,:\, \z_j>\vr\},
\]
where $\z_1\ge\ldots\ge\z_p$ are the positive eigenvalues of the Laplacian $\D_1$ and $p=\rank\D_1$.
\end{corollary}

\no \textbf{Proof.} Each fiber Laplacian $\D(\vt)$, $\vt\in\T^d$, is given by
\er{dec1}.  Then, due to \er{mpmm} -- \er{wnf1},
each eigenvalue $\l_j(\vt)$ of $\D(\vt)$ above its essential spectrum satisfy
\[\lb{com1}
\z_j\leq\z_j+m_-(\vt)\leq\l_j(\vt)\leq\z_j+m_+(\vt)\leq\z_j+\vr,
\]
which yields \er{esbp}. Let $\z_j>\vr$ for some $j=1,\ldots,p$.
Then, due to \er{com1},  $\l_j(\vt)>\vr$ for all $\vt\in\T^d$. Thus,
$\l_j$ creates the guided band $\gs_j^o(\D)=\l_j(\T^d)$. This yields
\er{esNg}. \qq \BBox

\

\no \textbf{Remark.} It is well known, see, e.g., \cite{F73}, that
the positive eigenvalues  $\z_1\ge\ldots\ge\z_p$ of the Laplacian
$\D_1$ on a finite graph $\G_1=(V_1,\cE_1)$ satisfy
$$
\frac{\nu_1}{\nu_1-1}\,\max_{v\in V_1}\vk^1_v\le\z_1\le \max_{u,v\in
V_1\atop u\sim v}(\vk^1_u+\vk^1_v),
$$
$$
\min_{u,v\in V_1\atop u\sim v}(\vk^1_u+\vk^1_v)-
(\nu_1-2)\le\z_p\le\frac{\nu_1}{\nu_1-1}\,\min_{v\in V_1}\vk^1_v,
$$
where $\nu_1=\# V_1$, $\vk_v^1$ is the degree of the vertex $v\in
V_1$ on $\G_1$.  From these estimates and \er{esbp} it follows that
increasing the degree of at least one vertex of the graph $\G_1$
removes the first guided band $\gs_1^o(\D)$ arbitrarily far to the
right. Increasing the degrees of all vertices of $\G_1$ removes all
guided bands arbitrarily far to the right.

\

%This Laplacian has the following decomposition:
%\[
%\lb{avo}
%\D^\gm=\wt\D^\gm_0+\wt\D_1,\qqq \wt\D^\gm_0=\D^\gm_0\oplus\O_{V_1\sm V^c_0},
%\qqq \wt\D_1=\O_{V^c_0\sm V_1}\oplus\D_1,
%\]
%where $\D^\gm_0$ is the Laplacian on the modified cylinder $\cC^\gm_0=(V_0^c,\cE_0^c\sm\cB^c)$.

\no \textbf{Proof of Theorem \ref{Est}.} i) We rewrite the fiber
Laplacian $\D(\vt)$,  $\vt\in\T^d$, defined by \er{l2.15''} in the
form:
\[
\label{eq.1}
\D(\vt)=\D^\gm+\D_\b(\vt),
\]
\[\label{Dbt}
\big(\D_\b(\vt)f\big)(v)=\sum_{\be=(v,\,u)\in\cA_c \atop \t(\be)\neq0}
\big(f(v)-e^{i\lan\t(\be),\,\vt\ran}f(u)\big), \qqq
 v\in V^c,
\]
where $\t(\be)\in\Z^d$ is the index of the edge $\be\in\cA^c$
defined by \er{in}, \er{inf}. Each operator $\D_\b(\vt)$,
$\vt\in\T^d$,  is the magnetic Laplacian on the graph
$\cC_\b=(V^c,\cB^c)$ and the degree of each vertex $v\in V^c$ on
$\cC_\b$ is equal to the number $\b_v$ of all bridges starting at
$v$. Then the  spectrum  $\s\big(\D_\b(\vt)\big)\ss[0,2\b_+]$ for each\, $\vt\in\T^d\ $
(see, e.g., \cite{HS99}) and,
due to \er{wnf1}, each eigenvalue $\l_j(\vt)$ of $\D(\vt)$ above its essential spectrum satisfy $\m_j\le\l_j(\vt)\le\m_j+2\b_+$, which yields \er{esbp1}.

ii) The existence of such graphs is proved in Propositions \ref{Prg1} and \ref{Prg3}. \qq \BBox

\subsection {Proof of Theorem \ref{T17}.}

Let $\G_t=(V_1,\cE_t)$ be a finite graph  obtained from the graph
$\G_1=(V_1,\cE_1)$ considering each edge of $\G_1$ to have the
multiplicity $t\in\N$. Then the Laplacian on the graph $\G_t$ has
the form $t\D_1$. If $t$ is large enough, then all positive
eigenvalues $t\z_p\le\ldots\le t\z_2\le t\z_1$ of the Laplacian
$t\D_1$ on the graph $\G_t$ satisfy $t\z_p>\vr$ and
$t(\z_j-\z_{j+1})>\vr$ for all $j\in\N_{p-1}$, where $\vr$ is
defined in \er{mm}. Then, due to Corollary \ref{Tloc}, the guided
spectrum of the Laplacian $\D_t$ consists of exactly $p$ guided
bands separated by gaps.

i) If $f_j=0$ on $V_{01}$, then, due to Proposition \ref{Pdis},
$\{t\z_j\}$ is a guided flat band of the Laplacian $\D_t$ on the
perturbed graph $\G=\G_0\cup\G_t^g$.

Let $f_j=(f_{j,01},f_{j,1})\in\ell^2(V_1)$, where $0\neq
f_{j,01}\in\ell^2(V_{01})$  and $f_{j,1}\in\ell^2(V_1\sm V_{01})$. Using \er{dec1}, we rewrite the fiber Laplacian $\D_t(\vt)$,
$\vt\in\T^d$,  for the Laplacian $\D_t$ acting on the perturbed
graph $\G=\G_0\cup\G_t^g$ in the form
$$
\D_t(\vt)=P\D_0(\vt)P+tP_1\D_1P_1=t K_t(\vt),\qq K_t(\vt)=P_1\D_1P_1+\ve P\D_0(\vt)P,\qq \ve=\frac1t\,.
$$
We denote the eigenvalues of the operator $K_t(\vt)$ above its
essential spectrum by
\[
\label{ebesH} E_p(\vt,t)\leq\ldots\leq E_2(\vt,t)\leq E_1(\vt,t),
\qqq \vt\in\T^d.
\]
Then the eigenvalue
$E_j(\vt,t)$ of the operator $K_t(\vt)$ has the following
asymptotics:
\[\lb{ass1}
E_j(\vt,t)=\z_j+\ve\,W_j(\vt)+O(\ve^2)
\]
(see pp. 7--8 in \cite{RS78}) uniformly in $\vt\in\T^d$ as $t\ra\iy$, where
\[
W_j(\vt)=\lan(0,f_{j,01}),\D_0(\vt)(0,f_{j,01})\ran_{V_0^c},
\]
$\lan\cdot\,,\cdot\ran_{V_0^c}$ denotes the inner product in
$\ell^2(V_0^c)$.  This yields the asymptotics of the eigenvalue
$\l_j(\vt,t)$ of the operator $\D_t(\vt)$:
\[\lb{ass}
\l_j(\vt,t)=t\,E_j(\vt,t)=t\z_j+W_j(\vt)+O(1/t).
\]
Using this asymptotics for
$\l_j^-(t)=\min\limits_{\vt\in\T^d}\l_j(\vt,t)$  and
$\l_j^+(t)=\max\limits_{\vt\in\T^d}\l_j(\vt,t)$, we obtain
$$
\l_j^\pm(t)=t\z_j+W_j^\pm+O(1/t),
$$
where $W_j^\pm$ are defined in \er{Dpm}. Since
$\gs_j(\D_t)=[\l_j^-(t),\l_j^+(t)]$, the asymptotics \er{ass} also
gives the  second formula in \er{Qt}. Using the formula \er{l2.15'}
for the fiber Laplacian $\D_0(\vt)$, we obtain
\[\lb{cjvt}
\begin{aligned}
&W_j(\vt)=\sum_{v\in V_{01}}\big(\D_0(\vt)(0,f_{j,01})\big)(v)\,\ol f_j(v)
=\sum_{v\in V_{01}}\sum_{\be=(v,\,u)\in\cA_0^c}
\big(f_j(v)-e^{i\lan\t(\be),\,\vt\ran}f_j(u)\big)\ol f_j(v)\\&=\sum_{v\in V_{01}}\vk_v^0\,|f_j(v)|^2-
\sum_{\be=(v,\,u)\in\cA_0^c \atop v,u\in V_{01}}e^{i\lan\t(\be),\,\vt\ran}f_j(u)\,\ol f_j(v)\\
&=\sum_{v\in V_{01}}\vk_v^0\,|f_j(v)|^2-
\sum_{\be=(v,\,u)\in\cA_0^c \atop v,u\in V_{01}}\cos\lan\t(\be),\vt\ran f_j(u)\,f_j(v),
\end{aligned}
\]
where $\vk_v^0$ is the degree of the vertex $v$ on the unperturbed
cylinder $\cC_0$  and $\t(\be)\in\Z^d$ is the edge index defined by
\er{in}, \er{inf}. Using this formula we rewrite the constant
$W_j^\bu$ defined in \er{Dpm} in the form
\[\lb{esCj}
W_j^\bu=\max_{\vt\in\T^d}\Omega_j(\vt)-\min_{\vt\in\T^d}\Omega_j(\vt),\qq
\Omega_j(\vt) =\sum_{\be=(v,\,u)\in\cA_0^c \atop v,u\in
V_{01},\t(\be)\neq0}\cos\lan\t(\be),\vt\ran f_j(u)\,f_j(v).
\]
We have
$$
|\Omega_j(\vt)|\leq\sum_{\be=(v,\,u)\in\cA_0^c \atop v,u\in V_{01},\t(\be)\neq0}|\cos\lan\t(\be),\vt
\ran|\leq\b_{01},
$$
which yields $W_j^\bu\leq2\b_{01}$.

ii) Let $V_{01}=\{v\}$. Then for $\Omega_j$ defined in \er{esCj} we have
\[\lb{V011}
\Omega_j(\vt)=f^2_j(v)\sum_{\be=(v,\,v)\in\cA_0^c \atop
\t(\be)\neq0}\cos\lan\t(\be),\vt\ran,  \qqq
\max_{\vt\in\T^d}\Omega_j(\vt)=f^2_j(v)\b_{01}.
\]
Using the identity $\int\limits_{\T^d}\cos\lan\t(\be),\vt\ran\,d\vt=0$ for each $\t(\be)\neq0$, we obtain
\[\lb{V012}
-f^2_j(v)\b_{01}\leq\min\limits_{\vt\in\T^d}\Omega_j(\vt)\leq0.
\]
Then \er{esCj} -- \er{V012} yield $f^2_j(v)\b_{01}\leq W_j^\bu\leq
2f^2_j(v)\b_{01}$.  From \er{V011} and the condition $f_j(v)\neq0$
it follows that $\Omega_j(\cdot)=\const$ iff $\b_{01}$=0. This
yields the last statement of the item.

iii) If all positive eigenvalues $\z_1>\ldots>\z_p$ of the
Laplacian $\D_1$ are distinct,  then summing the second asymptotics
in \er{Qt} over $j=1,\ldots,p$ we obtain \er{Qt1}.

iv) If there is no bridge connecting the vertices from the set
$V_{01}$ on  the cylinder $\cC$, then for each $j\in\N_p$ the
function $W_j$ defined by \er{cjvt} is constant, i.e., $W_j^\bu=0$,
and, the second asymptotics in \er{Qt} and the asymptotics \er{Qt1},
take the form $|\gs_j(\D_t)|=O(1/t)$ and $|\gs(\D_t)|=O(1/t)$,
respectively. \qq $\BBox$

\smallskip

\no \textbf{Remark.} The set $\cA_0^c$ of all oriented edges of the
cylinder $\cC_0$  is infinite, but the sum in \er{cjvt}  is taken
over a finite (maybe empty) set of edges $(v,u)\in\cA_0^c$ for which
$v,u\in V_{01}$.

\subsection{Geometric properties of the guided spectrum}
Now we prove Proposition \ref{Pdis} and Corollary \ref{TEg} about
geometric properties of the guided spectrum for specific graphs.

\smallskip

\no \textbf{Proof of Proposition \ref{Pdis}.} Let $\z$ be an
eigenvalue of  the Laplacian $\D_1$ on $\G_1$ with an eigenfunction
$f\in\ell^2(V_1)$ equal to zero on $V_{01}$. Then, due to \er{dec1},
for each $\vt\in\T^d$ the function $(0,f)\in\ell^2(V^c)$ is an
eigenfunction of the fiber Laplacian $\D(\vt)$ associated with the
same eigenvalue $\z$. Thus, $\{\z\}$ is a guided flat band of the
Laplacian $\D$. \qq $\BBox$

\smallskip

\no \textbf{Remarks.} 1) The $\n_1\ts\n_1$ matrix
$\D_1=\{\D^1_{uv}\}_{u,v\in V_1}$  associated to the Laplacian
$\D_1$ on a finite graph $\G_1=(V_1,\cE_1)$ in the standard
orthonormal basis is given by
\[\lb{MLsb}
\D^1_{uv}=\d_{uv}\vk_v^1-\vk_{uv}, \qqq \n_1=\# V_1,
\]
where $\d_{uv}$ is the Kronecker delta, $\vk_v^1$ is the degree of
the vertex $v\in V_1$  on the graph $\G_1$, $\vk_{uv}$ is the number
of edges $(u,v)\in\cA_1$.

2) The sufficient conditions from Proposition \ref{Pdis} are
equivalent to $\z$ being  an eigenvalue of two operators: the
Laplacian $\D_1$ and the Laplacian $\D_D$ on $\G_1$ with Dirichlet
boundary condition
\[
f(v)=0, \qqq \forall v\in V_{01},
\]
i.e., $\z$ is an eigenvalue of two matrices:
$\D_1=\{\D^1_{uv}\}_{u,v\in V_1}$ and its  submatrix
$\D_D=\{\D^1_{uv}\}_{u,v\in V_1\sm V_{01}}$.

\

\no \textbf{Proof of Corollary \ref{TEg}.} i) Due to the
connectivity of the periodic  graph $\G_0$, there exists a bridge on
the cylinder $\cC_0$. First we consider the case when there exists a
bridge-loop at some vertex $u_1\in V_0^c$. Due to the periodicity of
the cylinder $\cC_0$, for each $j\in\Z$ we have  $u_j=u_1+ja_{\wt
d}\in V_0^c$, where $a_{\wt d}$ is one of the periods of $\cC_0$.

\setlength{\unitlength}{1.0mm}
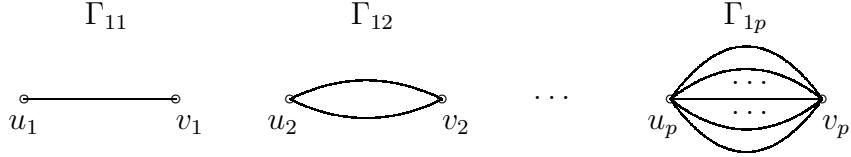
\begin{figure}[h]
\centering

\unitlength 1.0mm % = 2.845pt
\linethickness{0.4pt}
\ifx\plotpoint\undefined\newsavebox{\plotpoint}\fi % GNUPLOT compatibility
\begin{picture}(120,25)(0,0)

\put(5,10){\line(1,0){20.00}}
\put(5,10){\circle{1}}
\put(25,10){\circle{1}}
\put(13,20){$\G_{11}$}
\put(3,6){$u_1$}
\put(25,6){$v_1$}

\bezier{200}(40,10)(50,5)(60,10)
\bezier{200}(40,10)(50,15)(60,10)
\put(40,10){\circle{1}}
\put(60,10){\circle{1}}
\put(48,20){$\G_{12}$}
\put(37,6){$u_2$}
\put(60,6){$v_2$}

\put(72,10){$\ldots$}

\put(97,20){$\G_{1p}$}
\put(87,6){$u_p$}
\put(110,6){$v_p$}

\bezier{200}(90,10)(100,-4)(110,10)
\bezier{200}(90,10)(100,2)(110,10)
\put(90,10){\line(1,0){20.00}}
\bezier{200}(90,10)(100,18)(110,10)
\bezier{200}(90,10)(100,24)(110,10)
\put(90,10){\circle{1}}
\put(110,10){\circle{1}}
\put(98,8){$\ldots$}
\put(98,12){$\ldots$}
%***************************
\end{picture}

\vspace{-0.5cm} \caption{\footnotesize A finite graph $\G_1$ with
$p$ connected  components $\G_{11},\ldots,\G_{1p}$.} \label{gamt}
\end{figure}

Let $\G_1=(V_1,\cE_1)$ be a finite graph consisting of $p\in\N$
connected components $\G_{11},\ldots,\G_{1p}$, where $\G_{1j}$ is a
graph consisting of two vertices $u_j\in V_0^c$, $v_j\notin V_0^c$ and $j$ edges
connecting these vertices (see Fig.\ref{gamt}). The Laplacian
$\D_1$ on the finite graph $\G_1$ has exactly $p$ simple positive
eigenvalues
$$
\z_1=2p, \qq \z_2=2(p-1), \qq \ldots\,, \qq \z_p=2.
$$
The normalized eigenfunction $f_j\in\ell^2(V_1)$ of the Laplacian
$\D_1$ corresponding to the eigenvalue $\z_j$ has the form
\[\lb{egf}
f_j(v)=\left\{
\begin{array}{rl}
{1\/\sqrt{2}}\,, & \textrm{ if } \, v=u_j\\[4pt]
-{1\/\sqrt{2}}\,, & \textrm{ if } \, v=v_j\\[4pt]
0,  & \textrm{ otherwise}
\end{array}\right.,\qqq j=1,\ldots,p.
\]

Let $\G_t$, $t\in\N$, be a finite graph obtained from the graph
$\G_1$  considering each edge of $\G_1$ to have the multiplicity
$t$. Let $\D_t$ be the Laplacian on the perturbed graph
$\G=\G_0\cup\G_t^g$. Due to Theorem \ref{T17}, for $t$ large enough
the guided spectrum of the Laplacian $\D_t$ consists of exactly $p$
guided bands  $\gs_j^o(\D_t)=\gs_j(\D_t)$ separated by gaps and
these bands satisfy
\[\lb{gbnd}
|\gs_j(\D_t)|=W_j^\bu+O(1/t),\qqq j\in\N_p,
\]
where $W_j^\bu$ is defined in \er{Dpm}. Substituting the identities
\er{egf} into \er{cjvt}, we obtain the following expression for the
function $W_j$:
\[\lb{cjvt''}
W_j(\vt)={\vk_{u_j}^0\/2}-\sum_{\be=(u_j,\,u_j)\in\cA_0^c}\cos\lan\t(\be),\vt\ran.
\]
Due to the periodicity of the cylinder $\cC_0$, for each $\vt\in\T^d$ we have
\[\lb{D1p}
\sum_{\be=(u_1,\,u_1)\in\cA_0^c}\cos\lan\t(\be),\vt\ran=
\sum_{\be=(u_2,\,u_2)\in\cA_0^c}\cos\lan\t(\be),\vt\ran=\ldots\,.
\]
This and \er{cjvt''} yield that $W_1^\bu=\ldots=W_p^\bu\neq0$.
Then, due to \er{gbnd},  all guided bands are non-degenerate and
\[\lb{leme}
|\gs(\D_t)|=\sum\limits_{j=1}^p\big|\gs_j(\D_t)\big|=p\,W_1^\bu+O(1/t).
\]
Choosing $p>\dfrac C{W_1^\bu}$\,, we obtain $|\gs(\D_t)|>C$ for $t$ large enough.

Now let there be no bridge-loop on the cylinder $\cC_0$. Then the
existing   bridge connects some distinct vertices $u_1,v_1\in
V_0^c$. Repeating the above arguments but with $u_j=u_1+ja_{\wt
d}\in V_0^c$ and $v_j=v_1+ja_{\wt
d}\in V_0^c$, $j\in\N_p$, we also obtain the required statement.
Note that in this case the function $W_j$ has the form
\[
W_j(\vt)={1\/2}\big(\vk_{u_j}^0+\vk_{v_j}^0\big)+
\sum_{\be=(u_j,\,v_j)\in\cA_0^c}\cos\lan\t(\be),\vt\ran, \qqq j\in\N_p.
\]

ii) Let there exist a vertex $v\in V_0^c$ such that there is no
bridge on $\cC_0$ starting at $v$. We consider a finite graph
$\G_t=(V_1,\cE_t)$ consisting of two vertices $u\notin V_0^c$ and
$v$ and the edge $(u,v)$ of multiplicity $t\in\N$. The Laplacian on
the finite graph $\G_t$ has one simple positive eigenvalue
$\z_1=2t$. Let $\D_t$ be the Laplacian on the perturbed 
graph $\G=\G_0\cup\G_t^g$. Due to Theorem \ref{T17}, for $t$ large
enough the guided spectrum of $\D_t$ consists of exactly one guided
band $\gs_1(\D_t)$ and the length of this guided band satisfies
$|\gs_1(\D_t)|=O(1/t)$. Thus, for any small $\ve>0$ there exists
$t\in\N$ such that $|\gs(\D_t)|<\ve$.

iii) This follows from \er{esbp0} and the fact that the eigenvalues
of the Laplacian on the finite graph $\G_t$ can be arbitrary large as $t\ra\iy$.

iv) The sufficient conditions for the existence of   guided flat
bands are proved in Proposition \ref{Pdis}. For examples of finite
graphs $\G_1$ satisfying these conditions see Propositions
\ref{Prg1}, \ref{Prg3}. \qq \BBox

\subsection{Reduction to operators on unperturbed cylinder}
We reduce the eigenvalue problem for the fiber Laplacian on the
perturbed cylinder $\cC$ to that on the unperturbed cylinder
$\cC_0$. In order to do this we use the following well-known theorem \cite{BFS98}.

\begin{theorem}[Feshbach projection method]\lb{TFPM}
Let $P$ be an orthogonal projection on a separable Hilbert space
$\cH$, and  let $P^\perp=\1-P$ be its complement. Let $T$ be a
bounded self-adjoint operator. Assume that $P^\perp T P^\perp$ is
invertible on $P^\perp\cH$. Then
\begin{itemize}
  \item[i)] $T$ is invertible on $\cH$ if and only if its Feshbach map
\[\lb{fbm}
\cF=PTP-PTP^\perp\big(P^\perp TP^\perp\big)^{-1} P^\perp TP
\]
is invertible on $P\cH$; in this case $\cF^{-1}=PT^{-1}P$;
  \item[ii)] if $T\p=0$ for some vector $0\neq\p\in\cH$, then $\cF P\p=0$, where $P\p\neq0$;
  \item[iii)] if $\cF\vp=0$ for some vector $0\neq\vp=P\vp$, then $T\p=0$, where
      $$
      0\neq\p=\big[P- P^\perp(P^\perp T P^\perp)^{-1}P^\perp TP\big]\vp;
      $$
  \item[iv)] the kernels of $T$ and $\cF$ have equal dimensions.
\end{itemize}
\end{theorem}

\no \textbf{Remark.} The operator $T$ in our consideration is
$\D(\vt)-\l$,  where $\D(\vt)$ is the fiber Laplacian on the
perturbed cylinder $\cC=(V^c,\cE^c)$.

\begin{proposition}\label{PStSp}
Let $P$ be the orthogonal projection of $\ell^2(V^c)$ onto the
subspace $\ell^2(V_0^c)$.  Then the following statements hold true.

i) If $P^\perp(\D_1-\l)P^\perp$ is invertible on $P^\perp\ell^2(V^c)$, then
for each $\vt\in\T^d$
\[
\l\in\s_p\big(\D(\vt)\big)\qq \Leftrightarrow \qq 0\in\s_p\big(\cF(\vt,\l)\big),
\]
and the kernels of $\D(\vt)-\l$ and $\cF(\vt,\l)$ have equal
dimensions.  Here $\cF(\vt,\l)$ is the Feshbach map \er{fbm} for the
operator $\D(\vt)-\l$, defined by \er{l2.15''}, and $\cF(\vt,\l)$
has the form
\[\lb{sheM}
\cF(\vt,\l)=P\big(\D_0(\vt)-\l\big)P+P_{01}\big(\D_1-
\D_1P^\perp\big(P^\perp (\D_1-\l)P^\perp\big)^{-1}P^\perp\D_1\big)P_{01},
\]
$P_{01}$ is the orthogonal projection of $\ell^2(V^c)$ onto $\ell^2(V_{01})$.

ii) Let, in addition, a finite graph $\G_1$ consist of $c:=c_{\G_1}$ connected components
$$
\G_{11}=(V_{11},\cE_{11}),\qq \ldots\,, \qq \G_{1c}=(V_{1c},\cE_{1c})
$$
each of which has exactly one common vertex $v_1,\ldots,v_c$,
respectively,  with the unperturbed cylinder
$\cC_0=(V_0^c,\cE_0^c)$. Then for each $\vt\in\T^d$ the Feshbach map
\er{sheM} for the operator $\D(\vt)-\l$ is just a fiber
Schr\"odinger operator given by
\[\lb{sheSS}
\cF(\vt,\l)=P\big(\D_0(\vt)-\l\big)P+Q(\l),
\]
where $Q(\l)=Q(\l,\cdot)$ is a potential with the compact support $V_{01}=\{v_1,\ldots,v_c\}$:
\[\lb{cspQ}
Q(\l,v_j)=\big(\cP_j\big(\D_{1j}-\D_{1j}\cP_j^\perp\big(\cP_j^\perp
(\D_{1j}-\l)\cP_j^\perp\big)^{-1}\cP_j^\perp\D_{1j}\big)\cP_j\big)
\upharpoonright_{\ell^2(\{v_j\})},\qq j=1,\ldots,c,
\]
$\D_{1j}$ is the Laplacian on the finite graph $\G_{1j}$, $\cP_j$ is
the orthogonal projection of $\ell^2(V_{1j})$ onto the one-dimensional subspace $\ell^2(\{v_j\})$. In particular, if the operator $\D_{1j}-\l\cP_j^\perp$ is invertible, then the expression \er{cspQ} can
be written in the form
\[\lb{cspQ1}
Q^{-1}(\l,v_j)=\big(\cP_j(\D_{1j}-\l\cP_j^\perp)^{-1}\cP_j\big)
\upharpoonright_{\ell^2(\{v_j\})}.
\]

iii) If $P^\perp(\D_1-\l)P^\perp$ is not invertible on
$P^\perp\ell^2(V^c)$,  then $\l$ is a guided flat band of the
Laplacian $\D$ on the perturbed graph $\G=\G_0\cup\G_1^g$.
\end{proposition}

\no \textbf{Proof.} i) Let $P_1$ be the orthogonal projection of
$\ell^2(V^c)$  onto $\ell^2(V_1)$. Then, using \er{fbm}, \er{dec1}
and the identities
\[
PP_1=P_1P=P_{01},\qqq PP^\perp=P^\perp P=0, \qqq P^\perp P_1=P_1P^\perp=P^\perp,
\]
we have
\[
\begin{aligned}
&\cF(\vt,\l)=P\big(\D(\vt)-\l\big)P-P\D(\vt)P^\perp\big(P^\perp (\D(\vt)-\l)
P^\perp\big)^{-1}P^\perp\D(\vt)P
\\
&=P\big(P\D_0(\vt)P+P_1\D_1P_1-\l\big)P
\\
&-P\big(P\D_0(\vt)P+P_1\D_1P_1\big)P^\perp\big(P^\perp (P\D_0(\vt)P+
P_1\D_1P_1-\l)P^\perp\big)^{-1} P^\perp\big(P\D_0(\vt)P+P_1\D_1P_1\big)P
\\
&=P\big(\D_0(\vt)-\l\big)P+P_{01}\D_1P_{01}-P_{01}\D_1P^\perp\big(P^\perp
(\D_1-\l)P^\perp\big)^{-1}P^\perp\D_1P_{01}.
\end{aligned}
\]
This and Theorem \ref{TFPM} yield the required statement.

ii) In this case $\D_1=\D_{11}\oplus\ldots\oplus\D_{1c}$. Then
\er{sheM}  has the form \er{sheSS}, \er{cspQ}. Let the operator $\D_{1j}-\l\cP_j^\perp$ be invertible. Then, using the partitioned presentation of the inverse matrix (see p.18 in \cite{HJ85}), we obtain \er{cspQ1}.

iii) Let $P^\perp(\D_1-\l)P^\perp$ not be invertible on
$P^\perp\ell^2(V^c)$,  then $\l$ is an eigenvalue of the operator
$P^\perp\D_1P^\perp$ and, due to Proposition \ref{Pdis}, $\{\l\}$ is a guided
flat band of the Laplacian $\D$ on the perturbed graph
$\G$. \qq $\BBox$

\section{Square lattice with guides.}
\setcounter{equation}{0}
\lb{Sec5}

In this section we consider the square lattice with specific guides.
We obtain some properties of the guided spectrum on such graphs and
give  examples of the guided spectrum.

Let $\dL^2=(V,\cE)$ be the square lattice, where the vertex set
$V=\Z^2$ and the edge set $\cE=\big\{(m,m+(1,0)), (m,m+(0,1)),\;
\forall\,m\in\Z^2\big\}$,  see Fig.\ref{USL}.\emph{a}.
\setlength{\unitlength}{1.0mm}
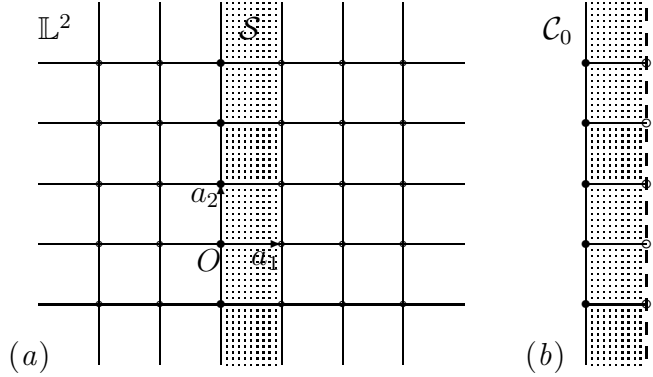
\begin{figure}[h]
\centering

\unitlength 0.8mm % = 2.845pt
\linethickness{0.4pt}
\ifx\plotpoint\undefined\newsavebox{\plotpoint}\fi % GNUPLOT compatibility
\begin{picture}(110,70)(0,0)

\put(-5,10){(\emph{a})}
\put(0,64){$\dL^2$}
\put(33,64){$\cS$}
%\put(77,38.5){$\Rightarrow$}

\put(0,20){\line(1,0){70.00}}
\put(0,30){\line(1,0){70.00}}
\put(0,40){\line(1,0){70.00}}
\put(0,50){\line(1,0){70.00}}
\put(0,60){\line(1,0){70.00}}

\put(10,10){\line(0,1){60.00}}
\put(20,10){\line(0,1){60.00}}
\put(30,10){\line(0,1){60.00}}
\put(40,10){\line(0,1){60.00}}
\put(50,10){\line(0,1){60.00}}
\put(60,10){\line(0,1){60.00}}

\bezier{60}(31,10)(31,40)(31,70)
\bezier{60}(32,10)(32,40)(32,70)
\bezier{60}(33,10)(33,40)(33,70)
\bezier{60}(34,10)(34,40)(34,70)
\bezier{60}(35,10)(35,40)(35,70)
\bezier{60}(36,10)(36,40)(36,70)
\bezier{60}(37,10)(37,40)(37,70)
\bezier{60}(38,10)(38,40)(38,70)
\bezier{60}(39,10)(39,40)(39,70)

\put(10,20){\circle{1}}
\put(20,20){\circle{1}}
\put(30,20){\circle*{1.5}}
\put(40,20){\circle{1}}
\put(50,20){\circle{1}}
\put(60,20){\circle{1}}

\put(10,30){\circle{1}}
\put(20,30){\circle{1}}
\put(30,30){\circle*{1.5}}
\put(40,30){\circle{1}}
\put(30,30){\vector(1,0){10.00}}
\put(30,30){\vector(0,1){10.00}}
\put(26,25.5){$O$}
\put(35,26.5){$a_1$}
\put(25.0,37.0){$a_2$}
\put(50,30){\circle{1}}
\put(60,30){\circle{1}}

\put(10,40){\circle{1}}
\put(20,40){\circle{1}}
\put(30,40){\circle*{1.5}}
\put(40,40){\circle{1}}
\put(50,40){\circle{1}}
\put(60,40){\circle{1}}

\put(10,50){\circle{1}}
\put(20,50){\circle{1}}
\put(30,50){\circle*{1.5}}
\put(40,50){\circle{1}}
\put(50,50){\circle{1}}
\put(60,50){\circle{1}}

\put(10,60){\circle{1}}
\put(20,60){\circle{1}}
\put(30,60){\circle*{1.5}}
\put(40,60){\circle{1}}
\put(50,60){\circle{1}}
\put(60,60){\circle{1}}

%***************************
\put(90,10){\line(0,1){60.0}}
\bezier{60}(91,10)(91,40)(91,70)
\bezier{60}(92,10)(92,40)(92,70)
\bezier{60}(93,10)(93,40)(93,70)
\bezier{60}(94,10)(94,40)(94,70)
\bezier{60}(95,10)(95,40)(95,70)
\bezier{60}(96,10)(96,40)(96,70)
\bezier{60}(97,10)(97,40)(97,70)
\bezier{60}(98,10)(98,40)(98,70)
\bezier{60}(99,10)(99,40)(99,70)
\multiput(100,11)(0,4){15}{\line(0,1){2}}

\put(90,20){\line(1,0){10.0}}
\put(90,30){\line(1,0){10.0}}
\put(90,40){\line(1,0){10.0}}
\put(90,50){\line(1,0){10.0}}
\put(90,60){\line(1,0){10.0}}

\put(90,20){\circle*{1.5}}
\put(100,20){\circle{1.5}}

\put(90,30){\circle*{1.5}}
\put(100,30){\circle{1.5}}

\put(90,40){\circle*{1.5}}
\put(100,40){\circle{1.5}}

\put(90,50){\circle*{1.5}}
\put(100,50){\circle{1.5}}

\put(90,60){\circle*{1.5}}
\put(100,60){\circle{1.5}}

\put(83,64){$\cC_0$}
\put(80,10){(\emph{b})}
%*******************************
\end{picture}
\vspace{-0.5cm} \caption{\footnotesize  \emph{a}) The square
lattice $\dL^2$; the vertices from the fundamental vertex set
$V_0^c$ are big black points; the strip $\cS$ is shaded;\; \emph{b})
the cylinder $\cC_0=\dL^2/\Z$ (the edges of the strip are
identified).}
\label{USL}
\end{figure}
The Laplacian $\D_0$ on $\dL^2$ has the form
\[\lb{Lao}
(\D_0 f)(m)=4f(m)-\sum_{|m-k|=1}f(k), \qqq f\in\ell^2(\Z^2), \qqq m\in\Z^2.
\]
We consider the Laplacian $\D$ on the perturbed square lattice
$\G=\dL^2\cup\G_1^g$ with  a guide $\G_1^g$. Due to Proposition
\ref{TSDI}, the Laplacian $\D$ on $\G$ has the decomposition
\er{raz}, \er{dec1} into a constant fiber direct integral, where the
fiber unperturbed Laplacian $\D_0(\vt)$ acts on $f\in\ell^2(\Z)$ and
is given by
\[
\label{l2.13}
\begin{aligned}
& \D_0(\vt)=2(1-\cos\vt)+h, \\
& (hf)(n)= 2f(n)-f(n+1)-f(n-1), \qqq n\in\Z,
\end{aligned}
\]
for all $\vt\in\T=(-\pi,\pi]$. It is well known that the spectrum of
the Laplacian $h$ on $\Z$ is given by $\s(h)=\s_{ac}(h)=[0,4]$. Then
 the spectrum of each fiber Laplacian $\D(\vt)$, $\vt\in\T$, has the form
\[\lb{sufl}
\begin{aligned}
&\s\big(\D(\vt)\big)=\s_{ac}\big(\D(\vt)\big)\cup\s_{p}\big(\D(\vt)\big),\\
& \s_{ac}\big(\D(\vt)\big)=\s_{ac}\big(\D_0(\vt)\big)
=[2-2\cos\vt,6-2\cos\vt],
\end{aligned}
\]
$\s_{p}\big(\D(\vt)\big)$ is the set of all eigenvalues of
$\D(\vt)$ of finite multiplicity given by \er{eq.3'}.

The spectrum of the Laplacian $\D$ on the perturbed square lattice
$\G=\dL^2\cup\G_1^g$ has the form
$$
\s(\D)=\s(\D_0)\cup\gs(\D), \qqq \s(\D_0)=[0,8],\qqq
\gs(\D)=\bigcup_{j=1}^{N}\gs_j(\D)=\gs_{ac}(\D)\cup\gs_{fb}(\D),
$$
where $N$ is defined in \er{sgSo}.

\medskip

Now we consider the perturbed square lattice $\G=\dL^2\cup\G_1^g$
when a finite graph $\G_1$ has exactly one common vertex $v_1=0$
with the unperturbed cylinder $\cC_0=\dL^2/\Z=(\Z,\cE_0^c)$, see Fig.~\ref{SLA}.

\setlength{\unitlength}{1.0mm}
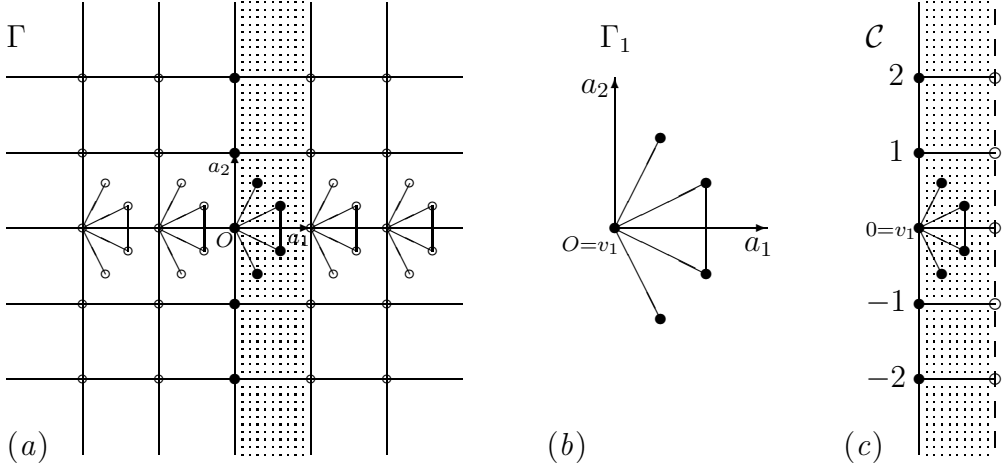
\begin{figure}[h]
\centering
\unitlength 1.0mm % = 2.845pt
\linethickness{0.4pt}
\ifx\plotpoint\undefined\newsavebox{\plotpoint}\fi % GNUPLOT compatibility
\begin{picture}(140,67)(0,0)
\put(0,10){(\emph{a})}
\put(0,64){$\G$}

\put(0,20){\line(1,0){60.00}}
\put(0,30){\line(1,0){60.00}}
\put(0,40){\line(1,0){60.00}}
\put(0,50){\line(1,0){60.00}}
\put(0,60){\line(1,0){60.00}}

\put(10,10){\line(0,1){60.00}}
\put(20,10){\line(0,1){60.00}}
\put(30,10){\line(0,1){60.00}}
\put(40,10){\line(0,1){60.00}}
\put(50,10){\line(0,1){60.00}}

\bezier{60}(31,10)(31,40)(31,70)
\bezier{60}(32,10)(32,40)(32,70)
\bezier{60}(33,10)(33,40)(33,70)
\bezier{60}(34,10)(34,40)(34,70)
\bezier{60}(35,10)(35,40)(35,70)
\bezier{60}(36,10)(36,40)(36,70)
\bezier{60}(37,10)(37,40)(37,70)
\bezier{60}(38,10)(38,40)(38,70)
\bezier{60}(39,10)(39,40)(39,70)

\put(10,20){\circle{1}}
\put(20,20){\circle{1}}
\put(30,20){\circle*{1.5}}
\put(40,20){\circle{1}}
\put(50,20){\circle{1}}

\put(10,30){\circle{1}}
\put(20,30){\circle{1}}
\put(30,30){\circle*{1.5}}
\put(40,30){\circle{1}}
\put(30,40){\vector(1,0){10.00}}
\put(30,40){\vector(0,1){10.00}}
\put(27.5,37.5){$\scriptstyle O$}
\put(37,38){$\scriptstyle a_1$}
\put(26.5,47.5){$\scriptstyle a_2$}
\put(50,30){\circle{1}}

\put(10,40){\circle{1}}
\put(20,40){\circle{1}}
\put(30,40){\circle*{1.5}}
\put(40,40){\circle{1}}
\put(50,40){\circle{1}}

\put(10,50){\circle{1}}
\put(20,50){\circle{1}}
\put(30,50){\circle*{1.5}}
\put(40,50){\circle{1}}
\put(50,50){\circle{1}}

\put(10,60){\circle{1}}
\put(20,60){\circle{1}}
\put(30,60){\circle*{1.5}}
\put(40,60){\circle{1}}
\put(50,60){\circle{1}}

\put(10,40){\line(1,2){3.00}}
\put(13,46){\circle{1}}
\put(10,40){\line(2,1){6.00}}
\put(16,43){\circle{1}}
\put(10,40){\line(1,-2){3.00}}
\put(13,34){\circle{1}}
\put(10,40){\line(2,-1){6.00}}
\put(16,37){\circle{1}}
\put(16,37){\line(0,1){6.00}}

\put(20,40){\line(1,2){3.00}}
\put(23,46){\circle{1}}
\put(20,40){\line(2,1){6.00}}
\put(26,43){\circle{1}}
\put(20,40){\line(1,-2){3.00}}
\put(23,34){\circle{1}}
\put(20,40){\line(2,-1){6.00}}
\put(26,37){\circle{1}}
\put(26,37){\line(0,1){6.00}}

\put(30,40){\line(1,2){3.00}}
\put(33,46){\circle*{1.5}}
\put(30,40){\line(2,1){6.00}}
\put(36,43){\circle*{1.5}}
\put(30,40){\line(1,-2){3.00}}
\put(33,34){\circle*{1.5}}
\put(30,40){\line(2,-1){6.00}}
\put(36,37){\circle*{1.5}}
\put(36,37){\line(0,1){6.00}}

\put(40,40){\line(1,2){3.00}}
\put(43,46){\circle{1}}
\put(40,40){\line(2,1){6.00}}
\put(46,43){\circle{1}}
\put(40,40){\line(1,-2){3.00}}
\put(43,34){\circle{1}}
\put(40,40){\line(2,-1){6.00}}
\put(46,37){\circle{1}}
\put(46,37){\line(0,1){6.00}}

\put(50,40){\line(1,2){3.00}}
\put(53,46){\circle{1}}
\put(50,40){\line(2,1){6.00}}
\put(56,43){\circle{1}}
\put(50,40){\line(1,-2){3.00}}
\put(53,34){\circle{1}}
\put(50,40){\line(2,-1){6.00}}
\put(56,37){\circle{1}}
\put(56,37){\line(0,1){6.00}}

%***********************
\put(80,40){\circle*{1.5}}
\put(80,40){\line(1,2){6.00}}
\put(86,52){\circle*{1.5}}
\put(80,40){\line(2,1){12.00}}
\put(92,46){\circle*{1.5}}
\put(80,40){\line(1,-2){6.00}}
\put(86,28){\circle*{1.5}}
\put(80,40){\line(2,-1){12.00}}
\put(92,34){\circle*{1.5}}
\put(92,34){\line(0,1){12.00}}
\put(80,40){\vector(1,0){20.00}}
\put(80,40){\vector(0,1){20.00}}
\put(73,37.0){$\scriptstyle O=v_1$}
\put(97,37){$a_1$}
\put(75.5,58.0){$a_2$}
\put(78.0,64){$\G_1$}
%\put(76.0,40){$v_1$}
\put(71.0,10){(\emph{b})}
%***************************
\put(120,10){\line(0,1){60.0}}
\multiput(130,11)(0,4){15}{\line(0,1){2}}

\bezier{60}(121,10)(121,40)(121,70)
\bezier{60}(122,10)(122,40)(122,70)
\bezier{60}(123,10)(123,40)(123,70)
\bezier{60}(124,10)(124,40)(124,70)
\bezier{60}(125,10)(125,40)(125,70)
\bezier{60}(126,10)(126,40)(126,70)
\bezier{60}(127,10)(127,40)(127,70)
\bezier{60}(128,10)(128,40)(128,70)
\bezier{60}(129,10)(129,40)(129,70)

\put(120,20){\line(1,0){10.0}}
\put(120,30){\line(1,0){10.0}}
\put(120,40){\line(1,0){10.0}}
\put(120,50){\line(1,0){10.0}}
\put(120,60){\line(1,0){10.0}}

\put(120,20){\circle*{1.5}}
\put(130,20){\circle{1.5}}

\put(120,30){\circle*{1.5}}
\put(130,30){\circle{1.5}}

\put(120,40){\circle*{1.5}}
\put(130,40){\circle{1.5}}

\put(120,50){\circle*{1.5}}
\put(130,50){\circle{1.5}}

\put(120,60){\circle*{1.5}}
\put(130,60){\circle{1.5}}

\put(120,40){\line(1,2){3.00}}
\put(123,46){\circle*{1.5}}
\put(120,40){\line(2,1){6.00}}
\put(126,43){\circle*{1.5}}
\put(120,40){\line(1,-2){3.00}}
\put(123,34){\circle*{1.5}}
\put(120,40){\line(2,-1){6.00}}
\put(126,37){\circle*{1.5}}
\put(126,37){\line(0,1){6.00}}
\put(113,19.0){$-2$}
\put(113,29.0){$-1$}
\put(113,39.0){$\scriptstyle 0=v_1$}
\put(116,49.0){$1$}
\put(116,59.0){$2$}
\put(113,64){$\cC$}
\put(110,10){(\emph{c})}
%*******************************
\end{picture}
\vspace{-0.5cm} \caption{\footnotesize  \emph{a}) The perturbed
square  lattice $\G=\dL^2\cup\G_1^g$; \; \emph{b}) a finite graph
$\G_1$; \; \emph{c})~the~perturbed cylinder $\cC=\G/\Z$.} \label{SLA}
\end{figure}

\begin{proposition}\label{Tnev1}
Let $\G_1=(V_1,\cE_1)$ be a finite connected graph and let
$V_{01}=\Z\cap V_1=\{0\}$. Then the  guided spectrum $\gs_{ac}(\D)$ and
$\gs_{fb}(\D)$ of the Laplacian $\D$ on the perturbed square lattice
$\G=\dL^2\cup\G_1^g$ satisfy
\[\lb{sgac}
\gs_{ac}(\D)=\Big\{\l\in\R\; : \; 2\leq\l-\sqrt{4+Q^2(\l)}\leq6,\qq
Q(\l)>0\Big\},
\]
\[\lb{sgfb}
\gs_{fb}(\D)=\big\{\l\in\R \; : \; \l \textrm{ is an
eigenvalue of the operator $\cP\D_1\cP$} \big\},
\]
where $\cP$ is the orthogonal projection of $\ell^2(V_1)$ onto the one-dimensional subspace $\ell^2(\{0\})$, $\D_1$ is the Laplacian on the finite graph $\G_1$,
\[\lb{ala}
Q(\l)=\big(\cP\big(\D_{1}-\D_{1}\cP^\perp\big(\cP^\perp (\D_{1}-\l)\cP^\perp\big)^{-1}\cP^\perp\D_{1}\big)\cP\big)
\upharpoonright_{\ell^2(\{0\})}.
\]
In particular, if the operator $\D_1-\l\cP^\perp$ is
invertible, then the expression \er{ala} can
be written in the form
\[\lb{ala1}
Q^{-1}(\l)=\big(\cP(\D_1-\l\cP^\perp)^{-1}\cP\big)
\upharpoonright_{\ell^2(\{0\})}.
\]
\end{proposition}

\no \textbf{Proof.} Let $\l\in\R$ not be an eigenvalue of the operator $\cP\D_1\cP$.
For each $\vt\in\T$,  using \er{sheSS}, \er{cspQ} and \er{l2.13}, we
obtain the Feshbach map for the operator $\D(\vt)-\l$:
\[\lb{heSS}
\begin{aligned}
& \cF(\vt,\l)=
P\big(\D_0(\vt)-\l\big)P+Q(\l)=P(h-\m)P+Q(\l),\\[6pt]
&\m=\l-(2-2\cos\vt),
\end{aligned}
\]
where $P$ is the orthogonal projection of $\ell^2(V^c)$ onto the
subspace $\ell^2(\Z)$,  $h$ is the Laplacian on $\Z$, the potential
$Q(\l)$ has the support $\{0\}$ and is given by \er{ala} or, in the case when $\D_1-\l\cP^\perp$ is
invertible, by \er{ala1}. Then, due to Proposition \ref{PStSp} and
the identity \er{heSS}, we have
\[\lb{coev}
\l\in\s_p\big(\D(\vt)\big)\qq \Leftrightarrow \qq
0\in\s_p\big(\cF(\vt,\l)\big) \qq \Leftrightarrow \qq
\m\in\s_p\big(h+Q(\l)\big).
\]
It is well known that the Schr\"odinger operator $h+Q$ on $\Z$ with
a potential $Q$ having  a support at a single point has the unique
eigenvalue
\[\lb{mu}
\m=\left\{
\begin{array}{rl}
2+\sqrt{4+Q^2}\;, & \qq Q>0, \\[6pt]
2-\sqrt{4+Q^2}\;, & \qq Q<0.
\end{array}\right.
\]
Combining this with \er{coev} and using the second identity in \er{heSS} and the fact that any eigenvalue $\l$ of $\D(\vt)$ satisfies $\l\geq2-2\cos\vt$, we obtain
\[\lb{sgac1}
\s_p\big(\D(\vt)\big)=\big\{\l\in\R : \; \l-\sqrt{4+Q^2(\l)}=4-2\cos\vt,\qq
Q(\l)>0\big\},
\]
which, due to the definition \er{dgsp} of the guided spectrum, yields \er{sgac} and $\l\notin\gs_{fb}(\D)$. 

Let $\l\in\R$ be an eigenvalue of the operator
$\cP\D_1\cP$.  Then, due to Proposition \ref{Pdis}, $\{\l\}$ is a guided
flat band of the Laplacian $\D$ on $\G$. Thus, $\gs_{fb}(\D)$ has the form \er{sgfb}.  \qq $\BBox$

\medskip

\no \textbf{Remark.} If $Q(\l)<0$ for some $\l>8$, then $\l$ lies in the gap of the spectrum  of the Laplacian $\D$ on the perturbed square lattice.

\subsection{Lattice with $p$ pendant edges at each vertex of $\Z\ts\{0\}$}
Let $\G_1=(V_1,\cE_1)$ be a star graph of order $p+1$, i.e., a tree
on $p+1$ vertices $v_1,v_2,\ldots,v_{p+1}$ with one vertex $v_1=0$
having degree $p$ and the other $p$ having degree 1. We assume that
\[\lb{stg}
V_{01}=\Z\cap V_1=\{0\}
\]
and consider the perturbed square lattice $\G=\dL^2\cup\G_1^g$, see Fig.\ref{SL}.

\setlength{\unitlength}{1.0mm}
\begin{figure}[h]
\centering

\unitlength 1.0mm % = 2.845pt
\linethickness{0.4pt}
\ifx\plotpoint\undefined\newsavebox{\plotpoint}\fi % GNUPLOT compatibility
\begin{picture}(130,70)(0,0)

\put(5,10){(\emph{a})}
\put(10,65){$\G$}
\put(10,20){\line(1,0){60.00}}
\put(10,30){\line(1,0){60.00}}
\put(10,40){\line(1,0){60.00}}
\put(10,50){\line(1,0){60.00}}
\put(10,60){\line(1,0){60.00}}

\put(20,10){\line(0,1){60.00}}
\put(30,10){\line(0,1){60.00}}
\put(40,10){\line(0,1){60.00}}
\put(50,10){\line(0,1){60.00}}
\put(60,10){\line(0,1){60.00}}

\bezier{60}(31,10)(31,40)(31,70)
\bezier{60}(32,10)(32,40)(32,70)
\bezier{60}(33,10)(33,40)(33,70)
\bezier{60}(34,10)(34,40)(34,70)
\bezier{60}(35,10)(35,40)(35,70)
\bezier{60}(36,10)(36,40)(36,70)
\bezier{60}(37,10)(37,40)(37,70)
\bezier{60}(38,10)(38,40)(38,70)
\bezier{60}(39,10)(39,40)(39,70)

\put(20,40){\line(1,2){3.00}}
\put(23,46){\circle{1}}
\put(20,40){\line(2,1){6.00}}
\put(26,43){\circle{1}}
\put(20,40){\line(1,-2){3.00}}
\put(23,34){\circle{1}}
\put(20,40){\line(2,-1){6.00}}
\put(26,37){\circle{1}}
\put(26,38){$\vdots$}

\put(30,40){\line(1,2){3.00}}
\put(33,46){\circle{1}}
\put(30,40){\line(2,1){6.00}}
\put(36,43){\circle{1}}
\put(30,40){\line(1,-2){3.00}}
\put(33,34){\circle{1}}
\put(30,40){\line(2,-1){6.00}}
\put(36,37){\circle{1}}
\put(36,38){$\vdots$}

\put(40,40){\line(1,2){3.00}}
\put(43,46){\circle{1}}
\put(40,40){\line(2,1){6.00}}
\put(46,43){\circle{1}}
\put(40,40){\line(1,-2){3.00}}
\put(43,34){\circle{1}}
\put(40,40){\line(2,-1){6.00}}
\put(46,37){\circle{1}}
\put(46,38){$\vdots$}

\put(50,40){\line(1,2){3.00}}
\put(53,46){\circle{1}}
\put(50,40){\line(2,1){6.00}}
\put(56,43){\circle{1}}
\put(50,40){\line(1,-2){3.00}}
\put(53,34){\circle{1}}
\put(50,40){\line(2,-1){6.00}}
\put(56,37){\circle{1}}
\put(56,38){$\vdots$}

\put(60,40){\line(1,2){3.00}}
\put(63,46){\circle{1}}
\put(60,40){\line(2,1){6.00}}
\put(66,43){\circle{1}}
\put(60,40){\line(1,-2){3.00}}
\put(63,34){\circle{1}}
\put(60,40){\line(2,-1){6.00}}
\put(66,37){\circle{1}}
\put(66,38){$\vdots$}

\put(20,20){\circle{1}}
\put(30,20){\circle{1}}
\put(40,20){\circle{1}}
\put(50,20){\circle{1}}
\put(60,20){\circle{1}}

\put(20,30){\circle{1}}
\put(30,30){\circle{1}}
\put(40,30){\circle{1}}
\put(50,30){\circle{1}}
\put(60,30){\circle{1}}

\put(20,40){\circle{1}}
\put(30,40){\circle{1}}
\put(40,40){\circle{1}}
\put(50,40){\circle{1}}
\put(60,40){\circle{1}}

\put(20,50){\circle{1}}
\put(30,50){\circle{1}}
\put(40,50){\circle{1}}
\put(50,50){\circle{1}}
\put(60,50){\circle{1}}

\put(20,60){\circle{1}}
\put(30,60){\circle{1}}
\put(40,60){\circle{1}}
\put(50,60){\circle{1}}
\put(60,60){\circle{1}}

%***************************

\put(100,40){\line(1,2){3.00}}
\put(103,46){\circle*{1.5}}
\put(100,40){\line(2,1){6.00}}
\put(106,43){\circle*{1.5}}
\put(100,40){\line(1,-2){3.00}}
\put(103,34){\circle*{1.5}}
\put(100,40){\line(2,-1){6.00}}
\put(106,37){\circle*{1.5}}
\put(106,38){$\vdots$}

\put(100,10){\line(0,1){60.0}}
\multiput(110,11)(0,4){15}{\line(0,1){2}}

\bezier{60}(101,10)(101,40)(101,70)
\bezier{60}(102,10)(102,40)(102,70)
\bezier{60}(103,10)(103,40)(103,70)
\bezier{60}(104,10)(104,40)(104,70)
\bezier{60}(105,10)(105,40)(105,70)
\bezier{60}(106,10)(106,40)(106,70)
\bezier{60}(107,10)(107,40)(107,70)
\bezier{60}(108,10)(108,40)(108,70)
\bezier{60}(109,10)(109,40)(109,70)

\put(100,20){\line(1,0){10.0}}
\put(100,30){\line(1,0){10.0}}
\put(100,40){\line(1,0){10.0}}
\put(100,50){\line(1,0){10.0}}
\put(100,60){\line(1,0){10.0}}

\put(100,20){\circle*{1.5}}
\put(110,20){\circle{1.5}}

\put(100,30){\circle*{1.5}}
\put(110,30){\circle{1.5}}

\put(100,40){\circle*{1.5}}
\put(110,40){\circle{1.5}}

\put(100,50){\circle*{1.5}}
\put(110,50){\circle{1.5}}

\put(100,60){\circle*{1.5}}
\put(110,60){\circle{1.5}}
\put(93,65){$\cC$}
\put(90,10){(\emph{b})}
\put(94,19.0){$-2$}
\put(94,29.0){$-1$}
\put(88,39.0){$0=v_1$}
%\put(97,39.0){$0$}
\put(97,49.0){$1$}
\put(97,59.0){$2$}
\put(103.5,47.0){$v_2$}
\put(106.5,43.5){$v_3$}
\put(106,34.0){$v_p$}
\put(102.5,31.0){$v_{p+1}$}
\end{picture}

\vspace{-0.5cm} \caption{\footnotesize  \emph{a}) The perturbed
square lattice $\G=\dL^2\cup\G_1^g$; \quad \emph{b}) the perturbed
cylinder $\cC=\G/\Z$.} \label{SL}
\end{figure}
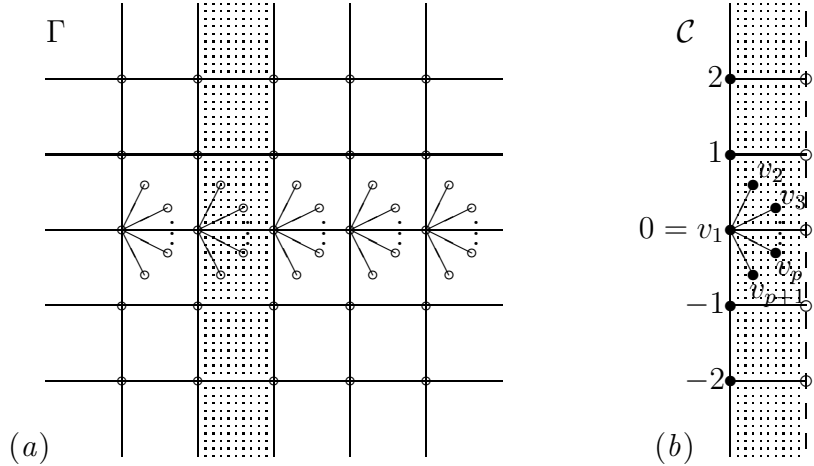

\begin{proposition}\label{Prg1}
The guided spectrum of the Laplacian $\D$ on the perturbed square
lattice $\G=\dL^2\cup\G_1^g$, where $\G_1$ is a star graph of order
$p+1$ satisfying \er{stg}, has the form
$$
\begin{aligned}
& \gs(\D)=\gs_{ac}(\D)\cup\gs_{fb}(\D),\qqq \gs_{ac}(\D)=\gs_1(\D),
\\
& \gs_{fb}(\D)=\gs_2(\D)\cup\ldots\cup\gs_{p}(\D)=\{1\},
\end{aligned}
$$
i.e., the guided flat band $\{1\}$ has the multiplicity $p-1$ and
$$
\gs_1(\D)\ss[p+3,p+8], \qqq |\gs_1(\D)|<4.
$$
Moreover, $\lim\limits_{p\ra\iy}|\gs_1(\D)|=2\b_+=4$, where  $\b_+$
is defined in \er{esbp1}.
\end{proposition}

\

\no\textbf{Remark.} For example, for $p=1$ and $p=2$ the  guided
band $\gs_1(\D)$ has the form
$$
\gs_1(\D)\approx[4{,}38;8{,}30] \qq \textrm{and} \qq
\gs_1(\D)\approx[5{,}18;9{,}01],
$$
respectively.

\no \textbf{Proof.} The $(p+1)\ts(p+1)$ matrix
$\D_1=\{\D^1_{uv}\}_{u,v\in V_1}$ defined by \er{MLsb} for the
Laplacian $\D_1$ on the star-graph $\G_1$ has the form
\[
\D_1=\left(
\begin{array}{cc}
  p & b\\
  b^* & I_p
\end{array}\right),\qqq b=(-1,\ldots,-1)\in\C^p,
\]
where $I_p$ is the identity $p\ts p$ matrix. Since $\l=1$ is an
eigenvalue  of the matrix $\D_1$ of multiplicity $p-1$, due to
Proposition \ref{Pdis}, the Laplacian $\D$ on $\G$ has a guided flat
band \{1\} of multiplicity $p-1$.

Now we find the absolutely continuous guided spectrum
$\gs_{ac}(\D)$.  Since the number of guided bands $N\leq p$ and the
guided flat band \{1\} has multiplicity $p-1$, the absolutely continuous
guided spectrum $\gs_{ac}(\D)$ consists of at most one band
$\gs_1(\D)=[\l_1^-,\l_1^+]$ and, due to Proposition \ref{Tnev1}, has the form \er{sgac}, where $Q(\l)={p\,\l\/\l-1}$\,. The inequality $Q(\l)>0$ gives that $\l>1$. The graphs of the function $f(\l)=\sqrt{4+Q^2(\l)}$ for $\l>1$ and of the functions $g_1(\l)=\l-2$, $g_2(\l)=\l-6$ are shown in Fig.\ref{Gfu}. Thus, the absolutely continuous guided spectrum $\gs_{ac}(\D)$ consists of exactly one guided band $\gs_1(\D)=[\l_1^-,\l_1^+]$ and $|\gs_1(\D)|=\l_1^+-\l_1^-<4$, see Fig.\ref{Gfu}.

\setlength{\unitlength}{1.0mm}
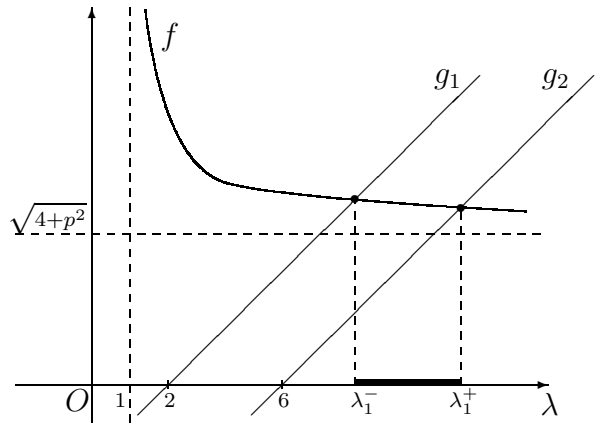
\begin{figure}[h]
\centering

\unitlength 1.0mm % = 2.845pt
\linethickness{0.4pt}
\ifx\plotpoint\undefined\newsavebox{\plotpoint}\fi % GNUPLOT compatibility
\begin{picture}(100,60)(0,0)

\put(0,10){\vector(1,0){70.00}}
\put(10,5){\vector(0,1){55.00}}
\put(20.0,9.5){\line(0,1){1}}
\put(35.0,9.5){\line(0,1){1}}
\put(13.0,7){$\scriptstyle 1$}
\put(19.5,7){$\scriptstyle 2$}
\put(34.5,7){$\scriptstyle 6$}
\multiput(15.0,5)(0,2){28}{\line(0,1){1}}
\multiput(0.0,30)(2,0){35}{\line(1,0){1}}

\put(44.6,34.6){\circle*{1}}
\put(58.5,33.5){\circle*{1}}
\multiput(44.6,10)(0,2){13}{\line(0,1){1}}
\multiput(58.5,10)(0,2){12}{\line(0,1){1}}
\put(44,7){$\scriptstyle \l_1^-$}
\put(57,7){$\scriptstyle \l_1^+$}
\put(44.6,10.2){\line(1,0){13.9}}
\put(44.6,10.4){\line(1,0){13.9}}
\put(44.6,10.6){\line(1,0){13.9}}

\put(6.5,6.5){$O$}
\put(-1.0,31.0){$\scriptstyle \sqrt{4+p^2}$}
\put(69,6){$\l$}
\put(55,50){$g_1$}
\put(69.2,50){$g_2$}
\put(19,55){$f$}

\put(16,6){\line(1,1){45.00}}
\put(31,6){\line(1,1){45.00}}

\bezier{700}(17.0,59.7)(19.5,40)(27.0,37)
\bezier{600}(27.0,37)(32,35)(67.0,33)

\end{picture}
\caption{\footnotesize  The graphs of the functions $f$, $g_1$, $g_2$.}
 \label{Gfu}
\end{figure}

From the identity $f(\l_1^-)=g_1(\l_1^-)$ we obtain that $\l_1^-$ is a root of the equation
\[\lb{eql0}
h(\l)\equiv\l^3-6\l^2+(9-p^2)\l-4=0.
\]
Using that
$$
h(p+3)=-4<0,\qqq h(p+4)=2p^2+9p>0,
$$
we obtain that the equation \er{eql0} has a unique solution
\[
\l_1^-\in[p+3,p+4].
\]
This and $|\gs_1(\D)|<4$ give that $\gs_1(\D)\ss[p+3,p+8]$.

We have $\l_1^\pm\ra\iy$ as $p\ra\iy$ and $f(\l)=\sqrt{4+p^2}+o(1)$
as $\l\ra\iy$. Then
$$
\l_1^-\ra2+\sqrt{4+p^2},\qqq \l_1^+\ra6+\sqrt{4+p^2},
$$
which yields that $|\gs_1(\D)|\ra4$ as $p\ra\iy$. On the other hand,
the number $\b_+=2$, since for each vertex of the cylinder $\cC$
there are two bridges starting at this vertex. Thus, we have
$|\gs_1(\D)|=2\b_++o(1)$ as $p\ra\iy$. \qq $\BBox$

\subsection{Lattice with mandarin guides} An \emph{$s$-mandarin graph}
 is a graph consisting of two vertices and $s$ edges connecting these vertices.
Let $\G_1=(V_1,\cE_1)$ be a union of two $s$-mandarin graphs with
one common vertex $v_1=0$, see Fig.\ref{SL2}.\emph{b}. We assume that
\[\lb{comv}
V_{01}=\Z\cap V_1=\{0\}
\]
and consider the perturbed square lattice $\G=\dL^2\cup\G_1^g$,  see
Fig.\ref{SL2}.

\begin{proposition}\label{Prg3}
The guided spectrum of the Laplacian $\D$ on the perturbed square
lattice $\G=\dL^2\cup\G_1^g$, where $\G_1$ is a union of two
$s$-mandarin graphs with one common vertex satisfying \er{comv}, has
the form
$$
\begin{aligned}
\gs(\D)=\gs_{ac}(\D)\cup\gs_{fb}(\D),\qqq \gs_{ac}(\D)=\gs_1(\D),
\qqq \gs_{fb}(\D)=\gs_2(\D)=\{s\},
\end{aligned}
$$
where
$$
\gs_1(\D)\ss[3s+1,3s+6]\qq \textrm{for} \qq s\geq2, \qqq |\gs_1(\D)|<4.
$$
Moreover, $\lim\limits_{s\ra\iy}|\gs_1(\D)|=4$.
\end{proposition}

\no\textbf{Remarks.} 1) For example, for $s=1$ and $s=2$ the guided
band  $\gs_1(\D)$ has the form
$$
\gs_1(\D)\approx[5{,}18;9{,}01] \qq \textrm{and} \qq
\gs_1(\D)\approx[7{,}75;11{,}26],
$$
respectively.

2) For $s\leq 8$ the guided flat band $\{s\}$ is embedded into the absolutely continuous spectrum of the Laplacian $\D$. For $s>8$ it lies in a gap.

\

\no \textbf{Proof.} The matrix $\D_1=\{\D^1_{uv}\}_{u,v\in V_1}$
defined by \er{MLsb} for the Laplacian $\D_1$ on the graph $\G_1$
has the form
\[
\D_1=s\left(
\begin{array}{cc}
  2 & b\\
  b^* & I_2
\end{array}\right),\qqq b=(-1,-1).
\]
Since $\l=s$ is a simple eigenvalue of the matrix $\D_1$, due to
Proposition \ref{Pdis}, the Laplacian $\D$ on $\G$ has a guided flat
band $\{s\}$.

The absolutely continuous guided spectrum $\gs_{ac}(\D)$ consists
of at most one band $\gs_1(\D)=[\l_1^-,\l_1^+]$ and, due to Proposition \ref{Tnev1}, has the form \er{sgac}, where $Q(\l)={2s\l\/\l-s}$\,. The inequality $Q(\l)>0$ gives that $\l>s$. The graphs of the function $f(\l)=\sqrt{4+Q^2(\l)}$ for $\l>s$ and of the functions $g_1(\l)=\l-2$, $g_2(\l)=\l-6$ are shown in Fig.\ref{Gfu2}. Thus, the absolutely
continuous guided spectrum $\gs_{ac}(\D)$ consists of exactly one
guided band $\gs_1(\D)=[\l_1^-,\l_1^+]$ and
$|\gs_1(\D)|=\l_1^+-\l_1^-<4$, see Fig.\ref{Gfu2}.

\setlength{\unitlength}{1.0mm}
\begin{figure}[h]
\centering

\unitlength 1.0mm % = 2.845pt
\linethickness{0.4pt}
\ifx\plotpoint\undefined\newsavebox{\plotpoint}\fi % GNUPLOT compatibility
\begin{picture}(100,60)(0,0)

\put(-10,10){\vector(1,0){105.00}}

\put(13,5){\vector(0,1){55.00}}
\put(20.0,9.5){\line(0,1){1}}
\put(35.0,9.5){\line(0,1){1}}

\put(19.5,7){$\scriptstyle 2$}
\put(34.5,7){$\scriptstyle 6$}
\put(41.0,7){$s$}
\multiput(40.0,5)(0,2){27}{\line(0,1){1}}
\multiput(-10.0,30)(2,0){52}{\line(1,0){1}}

\put(58.5,35.7){\circle*{1}}
\put(71.3,34.4){\circle*{1}}
\multiput(58.5,10)(0,2){13}{\line(0,1){1}}
\multiput(71.3,10.4)(0,2){12}{\line(0,1){1}}
\put(58,7){$\scriptstyle \l_1^-$}
\put(71,7){$\scriptstyle \l_1^+$}
\put(58.5,10.2){\line(1,0){12.8}}
\put(58.5,10.4){\line(1,0){12.8}}
\put(58.5,10.6){\line(1,0){12.8}}

\put(9.5,6.5){$O$}
\put(1.5,31.5){$\scriptstyle 2\sqrt{1+s^2}$}
\put(92,6){$\l$}
\put(75,47){$g_1$}
\put(90,47){$g_2$}
\put(44,55){$f$}
\put(13.8,6){\line(3,2){60}}
\put(28.8,6){\line(3,2){60}}

\bezier{700}(42.0,59.7)(44.5,40)(52.0,37)
\bezier{600}(52.0,37)(57,35)(92.0,33)

%\bezier{700}(38.0,59.7)(19.5,0)(8.0,17.0)
%\bezier{600}(8.0,17.0)(0,29)(-10.0,29)

\end{picture}
\caption{\footnotesize  The graphs of the functions $f$, $g_1$,
$g_2$.}  \label{Gfu2}
\end{figure}
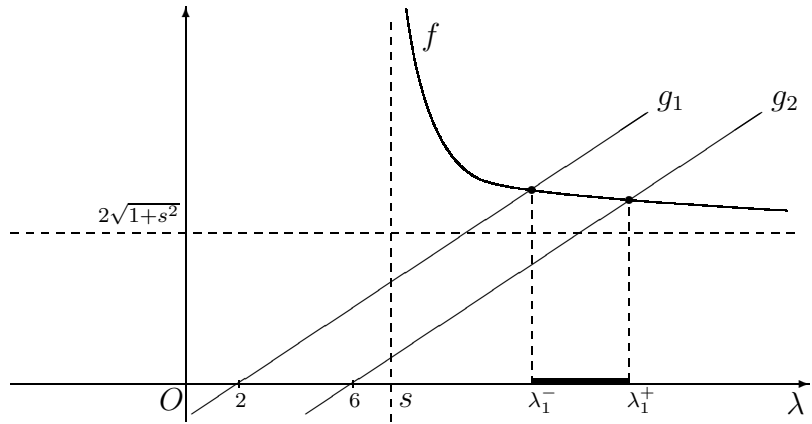

From the identity $f(\l_1^-)=g_1(\l_1^-)$ we obtain that $\l_1^-$ is a root of the equation
\[\lb{eql2}
h(\l)\equiv(\l-4)(\l-s)^2-4s^2\l=0.
\]
Using that
$$
h(3s+1)=-4s^2-9s-3<0,\qqq h(3s+2)=4(2s^2-s-2)>0 \qq \textrm{for} \qq s\geq2,
$$
we obtain that for $s\geq2$ the equation \er{eql2} has a unique solution
\[
\l_1^-\in[3s+1,3s+2].
\]
This and $|\gs_1(\D)|<4$ give that $\gs_1(\D)\ss[3s+1,3s+6]$.

We have $\l_1^\pm\ra\iy$ as $s\ra\iy$ and $f(\l)=2\sqrt{1+s^2}+o(1)$
as $\l\ra\iy$. Then
$$
\l_1^-\ra2+2\sqrt{1+s^2},\qqq \l_1^+\ra6+2\sqrt{1+s^2},
$$
which yields that $|\gs_1(\D)|\ra4$ as $s\ra\iy$. \qq $\BBox$

\subsection{Lattice with path guides} Let $\G_1=(V_1,\cE_1)$ be a path
of length 2, i.e., a connected graph with two vertices $v_1$ and
$v_3$ having degree 1 and a vertex $v_2$ having degree 2 and let
$\G_t=(V_1,\cE_t)$ be a finite graph obtained from the graph $\G_1$
considering each edge of $\G_1$ to have the multiplicity $t$. We
assume that $V_{01}=\Z\cap V_1=\{v_1=0\}$ and consider the perturbed square
lattice $\G=\dL^2\cup\G_t^g$, see Fig.\ref{SL3}.

\setlength{\unitlength}{1.0mm}
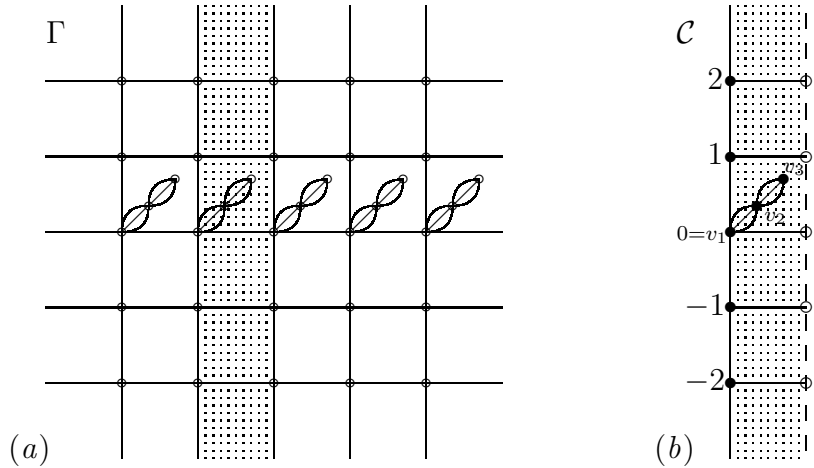
\begin{figure}[h]
\centering

\unitlength 1.0mm % = 2.845pt
\linethickness{0.4pt}
\ifx\plotpoint\undefined\newsavebox{\plotpoint}\fi % GNUPLOT compatibility
\begin{picture}(130,70)(0,0)

\put(5,10){(\emph{a})}
\put(10,65){$\G$}
\put(10,20){\line(1,0){60.00}}
\put(10,30){\line(1,0){60.00}}
\put(10,40){\line(1,0){60.00}}
\put(10,50){\line(1,0){60.00}}
\put(10,60){\line(1,0){60.00}}

\bezier{60}(31,10)(31,40)(31,70)
\bezier{60}(32,10)(32,40)(32,70)
\bezier{60}(33,10)(33,40)(33,70)
\bezier{60}(34,10)(34,40)(34,70)
\bezier{60}(35,10)(35,40)(35,70)
\bezier{60}(36,10)(36,40)(36,70)
\bezier{60}(37,10)(37,40)(37,70)
\bezier{60}(38,10)(38,40)(38,70)
\bezier{60}(39,10)(39,40)(39,70)

\put(20,10){\line(0,1){60.00}}
\put(30,10){\line(0,1){60.00}}
\put(40,10){\line(0,1){60.00}}
\put(50,10){\line(0,1){60.00}}
\put(60,10){\line(0,1){60.00}}

\put(20,40){\line(1,1){7.00}}
\put(23.5,43.5){\circle{1}}
\put(27,47){\circle{1}}
\bezier{100}(20,40)(20,43.5)(23.5,43.5)
\bezier{100}(20,40)(23.5,40)(23.5,43.5)
\bezier{100}(23.5,43.5)(23.5,47)(27,47)
\bezier{100}(23.5,43.5)(27,43.5)(27,47)

\put(30,40){\line(1,1){7.00}}
\put(33.5,43.5){\circle{1}}
\put(37,47){\circle{1}}
\bezier{100}(30,40)(30,43.5)(33.5,43.5)
\bezier{100}(30,40)(33.5,40)(33.5,43.5)
\bezier{100}(33.5,43.5)(33.5,47)(37,47)
\bezier{100}(33.5,43.5)(37,43.5)(37,47)

\put(40,40){\line(1,1){7.00}}
\put(43.5,43.5){\circle{1}}
\put(47,47){\circle{1}}
\bezier{100}(40,40)(40,43.5)(43.5,43.5)
\bezier{100}(40,40)(43.5,40)(43.5,43.5)
\bezier{100}(43.5,43.5)(43.5,47)(47,47)
\bezier{100}(43.5,43.5)(47,43.5)(47,47)

\put(50,40){\line(1,1){7.00}}
\put(53.5,43.5){\circle{1}}
\put(57,47){\circle{1}}
\bezier{100}(50,40)(50,43.5)(53.5,43.5)
\bezier{100}(50,40)(53.5,40)(53.5,43.5)
\bezier{100}(53.5,43.5)(53.5,47)(57,47)
\bezier{100}(53.5,43.5)(57,43.5)(57,47)

\put(60,40){\line(1,1){7.00}}
\put(63.5,43.5){\circle{1}}
\put(67,47){\circle{1}}
\bezier{100}(60,40)(60,43.5)(63.5,43.5)
\bezier{100}(60,40)(63.5,40)(63.5,43.5)
\bezier{100}(63.5,43.5)(63.5,47)(67,47)
\bezier{100}(63.5,43.5)(67,43.5)(67,47)

\put(20,20){\circle{1}}
\put(30,20){\circle{1}}
\put(40,20){\circle{1}}
\put(50,20){\circle{1}}
\put(60,20){\circle{1}}

\put(20,30){\circle{1}}
\put(30,30){\circle{1}}
\put(40,30){\circle{1}}
\put(50,30){\circle{1}}
\put(60,30){\circle{1}}

\put(20,40){\circle{1}}
\put(30,40){\circle{1}}
\put(40,40){\circle{1}}
\put(50,40){\circle{1}}
\put(60,40){\circle{1}}

\put(20,50){\circle{1}}
\put(30,50){\circle{1}}
\put(40,50){\circle{1}}
\put(50,50){\circle{1}}
\put(60,50){\circle{1}}

\put(20,60){\circle{1}}
\put(30,60){\circle{1}}
\put(40,60){\circle{1}}
\put(50,60){\circle{1}}
\put(60,60){\circle{1}}

%***************************

\put(100,10){\line(0,1){60.0}}
\multiput(110,11)(0,4){15}{\line(0,1){2}}

\bezier{60}(101,10)(101,40)(101,70)
\bezier{60}(102,10)(102,40)(102,70)
\bezier{60}(103,10)(103,40)(103,70)
\bezier{60}(104,10)(104,40)(104,70)
\bezier{60}(105,10)(105,40)(105,70)
\bezier{60}(106,10)(106,40)(106,70)
\bezier{60}(107,10)(107,40)(107,70)
\bezier{60}(108,10)(108,40)(108,70)
\bezier{60}(109,10)(109,40)(109,70)

\put(100,20){\line(1,0){10.0}}
\put(100,30){\line(1,0){10.0}}
\put(100,40){\line(1,0){10.0}}
\put(100,50){\line(1,0){10.0}}
\put(100,60){\line(1,0){10.0}}

\put(100,20){\circle*{1.5}}
\put(110,20){\circle{1.5}}

\put(100,30){\circle*{1.5}}
\put(110,30){\circle{1.5}}

\put(100,40){\circle*{1.5}}
\put(110,40){\circle{1.5}}

\put(100,50){\circle*{1.5}}
\put(110,50){\circle{1.5}}

\put(100,60){\circle*{1.5}}
\put(110,60){\circle{1.5}}
\put(93,65){$\cC$}
\put(90,10){(\emph{b})}
\put(94,19.0){$-2$}
\put(94,29.0){$-1$}
\put(93,39.0){$\scriptstyle 0=v_1$}
%\put(97,39.0){$0$}
\put(97,49.0){$1$}
\put(97,59.0){$2$}

\put(100,40){\line(1,1){7.00}}
\put(103.5,43.5){\circle*{1.5}}
\put(107,47){\circle*{1.5}}
\bezier{100}(100,40)(100,43.5)(103.5,43.5)
\bezier{100}(100,40)(103.5,40)(103.5,43.5)

\bezier{100}(103.5,43.5)(103.5,47)(107,47)
\bezier{100}(103.5,43.5)(107,43.5)(107,47)

\put(104.5,41.5){$\scriptstyle v_2$}
\put(107.0,48.0){$\scriptstyle v_3$}

\end{picture}

\vspace{-0.5cm} \caption{\footnotesize  \emph{a}) The perturbed
square lattice $\G=\dL^2\cup\G_t^g$; \quad \emph{b}) the perturbed
cylinder $\cC=\G/\Z$.} \label{SL3}
\end{figure}

\begin{proposition}\label{Prg2}
The guided spectrum of the Laplacian $\D$ on the perturbed square
lattice $\G=\dL^2\cup\G_t^g$ has the form
$$
\gs(\D)=\gs_{ac}(\D)=
\left\{
\begin{array}{ll}
 \gs_1(\D), & \textrm{ if } \; t=1,2\\[6pt]
 \gs_1(\D)\cup\gs_2(\D), & \textrm{ if } \; t\geq3
\end{array}\right.,
$$
where $\gs_2(\D)\ss[t,2t]$ and $\gs_1(\D)=[\l_1^-,\l_1^+]$,
\[\lb{ebb}
\l_1^-=2+\sqrt{4+t^2}+o(1),\qqq \l_1^+=6+\sqrt{4+t^2}+o(1),
\qq\textrm {as } \qq t\ra\iy.
\]
\end{proposition}

\no \textbf{Proof.} The matrix $\D_1=\{\D^1_{uv}\}_{u,v\in V_1}$
defined by \er{MLsb} for the Laplacian $\D_1$ on the graph $\G_t$
has the form
\[
\D_1=t\left(
\begin{array}{cc}
  1 & b \\
  b^* & \D_D
\end{array}\right),\qqq b=(-1,0),\qqq \D_D=\left(
\begin{array}{rr}
  2 & -1 \\
 -1 & 1
\end{array}\right).
\]
Since the matrices $\D_1$ and $\D_D$ have no equal eigenvalues, due
to Proposition \ref{Tnev1}, the Laplacian $\D$ on $\G$ has no guided
flat bands. Then the absolutely continuous guided spectrum
$\gs_{ac}(\D)$ consists of at most two bands and, due to Proposition \ref{Tnev1}, has the form \er{sgac}, where
$$
Q(\l)={t\l(\l-2t)\/\l^2-3t\l+t^2}\,.
$$
The inequality $Q(\l)>0$ gives that
\[\lb{a11p}
\textstyle \l\in\big({3-\sqrt{5}\/2}\,t,2t\big)
\cup\big({3+\sqrt{5}\/2}\,t,+\iy\big).
\]
For $t\geq4$ the graphs of the function $f=\sqrt{4+Q^2(\l)}$ for
$\l>{3-\sqrt{5}\/2}\,t$  and of the functions $g_1(\l)=\l-2$, $g_2(\l)=\l-6$ are shown in Fig.\ref{Gfu4}. Thus, for $t\geq4$ the absolutely
continuous guided spectrum $\gs_{ac}(\D)$ consists of exactly two
guided bands $\gs_s(\D)=[\l_s^-,\l_s^+]$, $s=1,2$. Since the
Laplacian $\D_1$ on the graph $\G_t$ has the positive eigenvalues
$\z_1=3t$ and $\z_2=t$, Corollary \ref{Tloc} and the formula
\er{a11p} give that $\gs_2(\D)\ss[t,2t]$.

We have $\l_1^\pm \ra\iy$ as $t\ra\iy$ and $f(\l)=\sqrt{4+t^2}+o(1)$
as $\l\ra\iy$. Then we get \er{ebb}.

From \er{sgac} by a direct calculation we obtain the absolutely
continuous guided spectrum of the Laplacian $\D$ for $t\leq 3$:
$$
\gs_{ac}(\D)=
\left\{
\begin{array}{ll}
 \gs_1(\D)\approx[4{,}47;8{,}31], \qqq & t=1\\[6pt]
 \gs_1(\D)\approx[6{,}66;9{,}45], \qqq & t=2\\[6pt]
 \gs_1(\D)\cup\gs_2(\D)\approx[4{,}62;6]\cup[9{,}51;11{,}44], \qqq & t=3\\[6pt]
\end{array}\right..
$$
$\BBox$

\setlength{\unitlength}{1.0mm}
\begin{figure}[h]
\centering

\unitlength 1.0mm % = 2.845pt
\linethickness{0.4pt}
\ifx\plotpoint\undefined\newsavebox{\plotpoint}\fi % GNUPLOT compatibility
\begin{picture}(100,60)(0,0)

\put(5,10){\vector(1,0){90.00}}

\put(13,5){\vector(0,1){55.00}}
\put(20.0,9.5){\line(0,1){1}}
\put(34.0,9.5){\line(0,1){1}}

\put(19.5,7){$\scriptstyle 2$}
\put(19.0,2.0){$\scriptstyle {3-\sqrt{5}\/2}\,t$}
\put(62.0,2.0){$\scriptstyle {3+\sqrt{5}\/2}\,t$}
\put(33.0,7){$\scriptstyle 6$}
\put(47.0,7){$\scriptstyle 2t$}
\multiput(48,10.5)(0,2){2}{\line(0,1){1}}
\multiput(22.0,5)(0,2){27}{\line(0,1){1}}
\multiput(65.0,5)(0,2){27}{\line(0,1){1}}
\multiput(5.0,30)(2,0){45}{\line(1,0){1}}

\put(35.7,17.7){\circle*{1}}
\put(43.3,14.5){\circle*{1}}
\multiput(35.7,10)(0,2){4}{\line(0,1){1}}
\multiput(43.3,10.0)(0,2){3}{\line(0,1){1}}
\put(35.5,7){$\scriptstyle \l_2^-$}
\put(42,7){$\scriptstyle \l_2^+$}
\put(35.7,10.2){\line(1,0){7.6}}
\put(35.7,10.4){\line(1,0){7.6}}
\put(35.7,10.6){\line(1,0){7.6}}

\put(75.8,37.6){\circle*{1}}
\put(84.3,34.9){\circle*{1}}
\multiput(75.8,10)(0,2){14}{\line(0,1){1}}
\multiput(84.3,10.4)(0,2){13}{\line(0,1){1}}
\put(75,7){$\scriptstyle \l_1^-$}
\put(84,7){$\scriptstyle \l_1^+$}
\put(75.8,10.2){\line(1,0){8.5}}
\put(75.8,10.4){\line(1,0){8.5}}
\put(75.8,10.6){\line(1,0){8.5}}

\put(9.5,6.5){$O$}
\put(2.5,32.0){$\scriptstyle \sqrt{4+t^2}$}
\put(92,6){$\l$}
\put(78,42){$g_1$}
\put(92,42){$g_2$}
\put(24,55){$f$}
\put(69,55){$f$}
\put(14.4,7){\line(2,1){67}}
\put(28.4,7){\line(2,1){67}}

\bezier{700}(67.0,58.0)(69.5,40)(77.0,37)
\bezier{600}(77.0,37)(82,35)(95.0,33)

\bezier{700}(63.5,58.0)(63.0,14)(48.0,14.0)
\bezier{600}(48.0,14.0)(23,14)(23.0,58)

\end{picture}
\caption{\footnotesize  The graphs of the functions $f$,  $g_1$,
$g_2$; $t\geq4$.} \label{Gfu4}
\end{figure}

%*****************************************************

\medskip

\footnotesize
 \textbf{Acknowledgments. \lb{Sec8}}  Our study was
supported by the RSF grant  No. 15-11-30007.

\end{document}